\input amstex 
%
\catcode`\@=11

\def\East#1#2{\setboxz@h{$\m@th\ssize\;{#1}\;\;$}%
 \setbox@ne\hbox{$\m@th\ssize\;{#2}\;\;$}\setbox\tw@\hbox{$\m@th#2$}%
 \dimen@\minaw@
 \ifdim\wdz@>\dimen@ \dimen@\wdz@ \fi  \ifdim\wd@ne>\dimen@ \dimen@\wd@ne \fi
 \ifdim\wd\tw@>\z@
  \mathrel{\mathop{\hbox to\dimen@{\rightarrowfill}}\limits^{#1}_{#2}}%
 \else
  \mathrel{\mathop{\hbox to\dimen@{\rightarrowfill}}\limits^{#1}}%
 \fi}
\def\West#1#2{\setboxz@h{$\m@th\ssize\;\;{#1}\;$}%
 \setbox@ne\hbox{$\m@th\ssize\;\;{#2}\;$}\setbox\tw@\hbox{$\m@th#2$}%
 \dimen@\minaw@
 \ifdim\wdz@>\dimen@ \dimen@\wdz@ \fi \ifdim\wd@ne>\dimen@ \dimen@\wd@ne \fi
 \ifdim\wd\tw@>\z@
  \mathrel{\mathop{\hbox to\dimen@{\leftarrowfill}}\limits^{#1}_{#2}}%
 \else
  \mathrel{\mathop{\hbox to\dimen@{\leftarrowfill}}\limits^{#1}}%
 \fi}
\font\arrow@i=lams1
\font\arrow@ii=lams2
\font\arrow@iii=lams3
\font\arrow@iv=lams4
\font\arrow@v=lams5
\newbox\zer@
\newdimen\standardcgap
\standardcgap=40\p@
\newdimen\hunit
\hunit=\tw@\p@
\newdimen\standardrgap
\standardrgap=32\p@
\newdimen\vunit
\vunit=1.6\p@
\def\Cgaps#1{\RIfM@
  \standardcgap=#1\standardcgap\relax \hunit=#1\hunit\relax
 \else \nonmatherr@\Cgaps \fi}
\def\Rgaps#1{\RIfM@
  \standardrgap=#1\standardrgap\relax \vunit=#1\vunit\relax
 \else \nonmatherr@\Rgaps \fi}
\newdimen\getdim@
\def\getcgap@#1{\ifcase#1\or\getdim@\z@\else\getdim@\standardcgap\fi}
\def\getrgap@#1{\ifcase#1\getdim@\z@\else\getdim@\standardrgap\fi}
\def\cgaps#1{\RIfM@
 \cgaps@{#1}\edef\getcgap@##1{\i@=##1\relax\the\toks@}\toks@{}\else
 \nonmatherr@\cgaps\fi}
\def\rgaps#1{\RIfM@
 \rgaps@{#1}\edef\getrgap@##1{\i@=##1\relax\the\toks@}\toks@{}\else
 \nonmatherr@\rgaps\fi}
\def\Gaps@@{\gaps@@}
\def\cgaps@#1{\toks@{\ifcase\i@\or\getdim@=\z@}%
 \gaps@@\standardcgap#1;\gaps@@\gaps@@
 \edef\next@{\the\toks@\noexpand\else\noexpand\getdim@\noexpand\standardcgap
  \noexpand\fi}%
 \toks@=\expandafter{\next@}}
\def\rgaps@#1{\toks@{\ifcase\i@\getdim@=\z@}%
 \gaps@@\standardrgap#1;\gaps@@\gaps@@
 \edef\next@{\the\toks@\noexpand\else\noexpand\getdim@\noexpand\standardrgap
  \noexpand\fi}%
 \toks@=\expandafter{\next@}}
\def\gaps@@#1#2;#3{\mgaps@#1#2\mgaps@
 \edef\next@{\the\toks@\noexpand\or\noexpand\getdim@
  \noexpand#1\the\mgapstoks@@}%
 \global\toks@=\expandafter{\next@}%
 \DN@{#3}%
 \ifx\next@\Gaps@@\gdef\next@##1\gaps@@{}\else
  \gdef\next@{\gaps@@#1#3}\fi\next@}
\def\mgaps@#1{\let\mgapsnext@#1\FN@\mgaps@@}
\def\mgaps@@{\ifx\next\space@\DN@. {\FN@\mgaps@@}\else
 \DN@.{\FN@\mgaps@@@}\fi\next@.}
\def\mgaps@@@{\ifx\next\w\let\next@\mgaps@@@@\else
 \let\next@\mgaps@@@@@\fi\next@}
\newtoks\mgapstoks@@
\def\mgaps@@@@@#1\mgaps@{\getdim@\mgapsnext@\getdim@#1\getdim@
 \edef\next@{\noexpand\getdim@\the\getdim@}%
 \mgapstoks@@=\expandafter{\next@}}
\def\mgaps@@@@\w#1#2\mgaps@{\mgaps@@@@@#2\mgaps@
 \setbox\zer@\hbox{$\m@th\hskip15\p@\tsize@#1$}%
 \dimen@\wd\zer@
 \ifdim\dimen@>\getdim@ \getdim@\dimen@ \fi
 \edef\next@{\noexpand\getdim@\the\getdim@}%
 \mgapstoks@@=\expandafter{\next@}}
\def\changewidth#1#2{\setbox\zer@\hbox{$\m@th#2}%
 \hbox to\wd\zer@{\hss$\m@th#1$\hss}}
\atdef@({\FN@\ARROW@}
\def\ARROW@{\ifx\next)\let\next@\OPTIONS@\else
 \DN@{\csname\string @(\endcsname}\fi\next@}
\newif\ifoptions@
\def\OPTIONS@){\ifoptions@\let\next@\relax\else
 \DN@{\options@true\begingroup\optioncodes@}\fi\next@}
\newif\ifN@
\newif\ifE@
\newif\ifNESW@
\newif\ifH@
\newif\ifV@
\newif\ifHshort@
\expandafter\def\csname\string @(\endcsname #1,#2){%
 \ifoptions@\let\next@\endgroup\else\let\next@\relax\fi\next@
 \N@false\E@false\H@false\V@false\Hshort@false
 \ifnum#1>\z@\E@true\fi
 \ifnum#1=\z@\V@true\tX@false\tY@false\a@false\fi
 \ifnum#2>\z@\N@true\fi
 \ifnum#2=\z@\H@true\tX@false\tY@false\a@false\ifshort@\Hshort@true\fi\fi
 \NESW@false
 \ifN@\ifE@\NESW@true\fi\else\ifE@\else\NESW@true\fi\fi
 \arrow@{#1}{#2}%
 \global\options@false
 \global\scount@\z@\global\tcount@\z@\global\arrcount@\z@
 \global\s@false\global\sxdimen@\z@\global\sydimen@\z@
 \global\tX@false\global\tXdimen@i\z@\global\tXdimen@ii\z@
 \global\tY@false\global\tYdimen@i\z@\global\tYdimen@ii\z@
 \global\a@false\global\exacount@\z@
 \global\x@false\global\xdimen@\z@
 \global\X@false\global\Xdimen@\z@
 \global\y@false\global\ydimen@\z@
 \global\Y@false\global\Ydimen@\z@
 \global\p@false\global\pdimen@\z@
 \global\label@ifalse\global\label@iifalse
 \global\dl@ifalse\global\ldimen@i\z@
 \global\dl@iifalse\global\ldimen@ii\z@
 \global\short@false\global\unshort@false}
\newif\iflabel@i
\newif\iflabel@ii
\newcount\scount@
\newcount\tcount@
\newcount\arrcount@
\newif\ifs@
\newdimen\sxdimen@
\newdimen\sydimen@
\newif\iftX@
\newdimen\tXdimen@i
\newdimen\tXdimen@ii
\newif\iftY@
\newdimen\tYdimen@i
\newdimen\tYdimen@ii
\newif\ifa@
\newcount\exacount@
\newif\ifx@
\newdimen\xdimen@
\newif\ifX@
\newdimen\Xdimen@
\newif\ify@
\newdimen\ydimen@
\newif\ifY@
\newdimen\Ydimen@
\newif\ifp@
\newdimen\pdimen@
\newif\ifdl@i
\newif\ifdl@ii
\newdimen\ldimen@i
\newdimen\ldimen@ii
\newif\ifshort@
\newif\ifunshort@
\def\zero@#1{\ifnum\scount@=\z@
 \if#1e\global\scount@\m@ne\else
 \if#1t\global\scount@\tw@\else
 \if#1h\global\scount@\thr@@\else
 \if#1'\global\scount@6 \else
 \if#1`\global\scount@7 \else
 \if#1(\global\scount@8 \else
 \if#1)\global\scount@9 \else
 \if#1s\global\scount@12 \else
 \if#1H\global\scount@13 \else
 \Err@{\Invalid@@ option \string\0}\fi\fi\fi\fi\fi\fi\fi\fi\fi
 \fi}
\def\one@#1{\ifnum\tcount@=\z@
 \if#1e\global\tcount@\m@ne\else
 \if#1h\global\tcount@\tw@\else
 \if#1t\global\tcount@\thr@@\else
 \if#1'\global\tcount@4 \else
 \if#1`\global\tcount@5 \else
 \if#1(\global\tcount@10 \else
 \if#1)\global\tcount@11 \else
 \if#1s\global\tcount@12 \else
 \if#1H\global\tcount@13 \else
 \Err@{\Invalid@@ option \string\1}\fi\fi\fi\fi\fi\fi\fi\fi\fi
 \fi}
\def\a@#1{\ifnum\arrcount@=\z@
 \if#10\global\arrcount@\m@ne\else
 \if#1+\global\arrcount@\@ne\else
 \if#1-\global\arrcount@\tw@\else
 \if#1=\global\arrcount@\thr@@\else
 \Err@{\Invalid@@ option \string\a}\fi\fi\fi\fi
 \fi}
\def\ds@(#1;#2){\ifs@\else
 \global\s@true
 \sxdimen@\hunit \global\sxdimen@#1\sxdimen@\relax
 \sydimen@\vunit \global\sydimen@#2\sydimen@\relax
 \fi}
\def\dtX@(#1;#2){\iftX@\else
 \global\tX@true
 \tXdimen@i\hunit \global\tXdimen@i#1\tXdimen@i\relax
 \tXdimen@ii\vunit \global\tXdimen@ii#2\tXdimen@ii\relax
 \fi}
\def\dtY@(#1;#2){\iftY@\else
 \global\tY@true
 \tYdimen@i\hunit \global\tYdimen@i#1\tYdimen@i\relax
 \tYdimen@ii\vunit \global\tYdimen@ii#2\tYdimen@ii\relax
 \fi}
\def\da@#1{\ifa@\else\global\a@true\global\exacount@#1\relax\fi}
\def\dx@#1{\ifx@\else
 \global\x@true
 \xdimen@\hunit \global\xdimen@#1\xdimen@\relax
 \fi}
\def\dX@#1{\ifX@\else
 \global\X@true
 \Xdimen@\hunit \global\Xdimen@#1\Xdimen@\relax
 \fi}
\def\dy@#1{\ify@\else
 \global\y@true
 \ydimen@\vunit \global\ydimen@#1\ydimen@\relax
 \fi}
\def\dY@#1{\ifY@\else
 \global\Y@true
 \Ydimen@\vunit \global\Ydimen@#1\Ydimen@\relax
 \fi}
\def\p@@#1{\ifp@\else
 \global\p@true
 \pdimen@\hunit \divide\pdimen@\tw@ \global\pdimen@#1\pdimen@\relax
 \fi}
\def\L@#1{\iflabel@i\else
 \global\label@itrue \gdef\label@i{#1}%
 \fi}
\def\l@#1{\iflabel@ii\else
 \global\label@iitrue \gdef\label@ii{#1}%
 \fi}
\def\dL@#1{\ifdl@i\else
 \global\dl@itrue \ldimen@i\hunit \global\ldimen@i#1\ldimen@i\relax
 \fi}
\def\dl@#1{\ifdl@ii\else
 \global\dl@iitrue \ldimen@ii\hunit \global\ldimen@ii#1\ldimen@ii\relax
 \fi}
\def\s@{\ifunshort@\else\global\short@true\fi}
\def\uns@{\ifshort@\else\global\unshort@true\global\short@false\fi}
\def\optioncodes@{\let\0\zero@\let\1\one@\let\a\a@\let\ds\ds@\let\dtX\dtX@
 \let\dtY\dtY@\let\da\da@\let\dx\dx@\let\dX\dX@\let\dY\dY@\let\dy\dy@
 \let\p\p@@\let\L\L@\let\l\l@\let\dL\dL@\let\dl\dl@\let\s\s@\let\uns\uns@}
\def\slopes@{\\161\\152\\143\\134\\255\\126\\357\\238\\349\\45{10}\\56{11}%
 \\11{12}\\65{13}\\54{14}\\43{15}\\32{16}\\53{17}\\21{18}\\52{19}\\31{20}%
 \\41{21}\\51{22}\\61{23}}
\newcount\tan@i
\newcount\tan@ip
\newcount\tan@ii
\newcount\tan@iip
\newdimen\slope@i
\newdimen\slope@ip
\newdimen\slope@ii
\newdimen\slope@iip
\newcount\angcount@
\newcount\extracount@
\def\slope@{{\slope@i=\secondy@ \advance\slope@i-\firsty@
 \ifN@\else\multiply\slope@i\m@ne\fi
 \slope@ii=\secondx@ \advance\slope@ii-\firstx@
 \ifE@\else\multiply\slope@ii\m@ne\fi
 \ifdim\slope@ii<\z@
  \global\tan@i6 \global\tan@ii\@ne \global\angcount@23
 \else
  \dimen@\slope@i \multiply\dimen@6
  \ifdim\dimen@<\slope@ii
   \global\tan@i\@ne \global\tan@ii6 \global\angcount@\@ne
  \else
   \dimen@\slope@ii \multiply\dimen@6
   \ifdim\dimen@<\slope@i
    \global\tan@i6 \global\tan@ii\@ne \global\angcount@23
   \else
    \tan@ip\z@ \tan@iip \@ne
    \def\\##1##2##3{\global\angcount@=##3\relax
     \slope@ip\slope@i \slope@iip\slope@ii
     \multiply\slope@iip##1\relax \multiply\slope@ip##2\relax
     \ifdim\slope@iip<\slope@ip
      \global\tan@ip=##1\relax \global\tan@iip=##2\relax
     \else
      \global\tan@i=##1\relax \global\tan@ii=##2\relax
      \def\\####1####2####3{}%
     \fi}%
    \slopes@
    \slope@i=\secondy@ \advance\slope@i-\firsty@
    \ifN@\else\multiply\slope@i\m@ne\fi
    \multiply\slope@i\tan@ii \multiply\slope@i\tan@iip \multiply\slope@i\tw@
    \count@\tan@i \multiply\count@\tan@iip
    \extracount@\tan@ip \multiply\extracount@\tan@ii
    \advance\count@\extracount@
    \slope@ii=\secondx@ \advance\slope@ii-\firstx@
    \ifE@\else\multiply\slope@ii\m@ne\fi
    \multiply\slope@ii\count@
    \ifdim\slope@i<\slope@ii
     \global\tan@i=\tan@ip \global\tan@ii=\tan@iip
     \global\advance\angcount@\m@ne
    \fi
   \fi
  \fi
 \fi}%
}
\def\slope@a#1{{\def\\##1##2##3{\ifnum##3=#1\global\tan@i=##1\relax
 \global\tan@ii=##2\relax\fi}\slopes@}}
\newcount\i@
\newcount\j@
\newcount\colcount@
\newcount\Colcount@
\newcount\tcolcount@
\newdimen\rowht@
\newdimen\rowdp@
\newcount\rowcount@
\newcount\Rowcount@
\newcount\maxcolrow@
\newtoks\colwidthtoks@
\newtoks\Rowheighttoks@
\newtoks\Rowdepthtoks@
\newtoks\widthtoks@
\newtoks\Widthtoks@
\newtoks\heighttoks@
\newtoks\Heighttoks@
\newtoks\depthtoks@
\newtoks\Depthtoks@
\newif\iffirstnewCDcr@
\def\dotoks@i{%
 \global\widthtoks@=\expandafter{\the\widthtoks@\else\getdim@\z@\fi}%
 \global\heighttoks@=\expandafter{\the\heighttoks@\else\getdim@\z@\fi}%
 \global\depthtoks@=\expandafter{\the\depthtoks@\else\getdim@\z@\fi}}
\def\dotoks@ii{%
 \global\widthtoks@{\ifcase\j@}%
 \global\heighttoks@{\ifcase\j@}%
 \global\depthtoks@{\ifcase\j@}}
\def\prenewCD@#1\endnewCD{\setbox\zer@
 \vbox{%
  \def\arrow@##1##2{{}}%
  \rowcount@\m@ne \colcount@\z@ \Colcount@\z@
  \firstnewCDcr@true \toks@{}%
  \widthtoks@{\ifcase\j@}%
  \Widthtoks@{\ifcase\i@}%
  \heighttoks@{\ifcase\j@}%
  \Heighttoks@{\ifcase\i@}%
  \depthtoks@{\ifcase\j@}%
  \Depthtoks@{\ifcase\i@}%
  \Rowheighttoks@{\ifcase\i@}%
  \Rowdepthtoks@{\ifcase\i@}%
  \Let@
  \everycr{%
   \noalign{%
    \global\advance\rowcount@\@ne
    \ifnum\colcount@<\Colcount@
    \else
     \global\Colcount@=\colcount@ \global\maxcolrow@=\rowcount@
    \fi
    \global\colcount@\z@
    \iffirstnewCDcr@
     \global\firstnewCDcr@false
    \else
     \edef\next@{\the\Rowheighttoks@\noexpand\or\noexpand\getdim@\the\rowht@}%
      \global\Rowheighttoks@=\expandafter{\next@}%
     \edef\next@{\the\Rowdepthtoks@\noexpand\or\noexpand\getdim@\the\rowdp@}%
      \global\Rowdepthtoks@=\expandafter{\next@}%
     \global\rowht@\z@ \global\rowdp@\z@
     \dotoks@i
     \edef\next@{\the\Widthtoks@\noexpand\or\the\widthtoks@}%
      \global\Widthtoks@=\expandafter{\next@}%
     \edef\next@{\the\Heighttoks@\noexpand\or\the\heighttoks@}%
      \global\Heighttoks@=\expandafter{\next@}%
     \edef\next@{\the\Depthtoks@\noexpand\or\the\depthtoks@}%
      \global\Depthtoks@=\expandafter{\next@}%
     \dotoks@ii
    \fi}}%
  \tabskip\z@
  \halign{&\setbox\zer@\hbox{\vrule height10\p@ width\z@ depth\z@
   $\m@th\displaystyle{##}$}\copy\zer@
   \ifdim\ht\zer@>\rowht@ \global\rowht@\ht\zer@ \fi
   \ifdim\dp\zer@>\rowdp@ \global\rowdp@\dp\zer@ \fi
   \global\advance\colcount@\@ne
   \edef\next@{\the\widthtoks@\noexpand\or\noexpand\getdim@\the\wd\zer@}%
    \global\widthtoks@=\expandafter{\next@}%
   \edef\next@{\the\heighttoks@\noexpand\or\noexpand\getdim@\the\ht\zer@}%
    \global\heighttoks@=\expandafter{\next@}%
   \edef\next@{\the\depthtoks@\noexpand\or\noexpand\getdim@\the\dp\zer@}%
    \global\depthtoks@=\expandafter{\next@}%
   \cr#1\crcr}}%
 \Rowcount@=\rowcount@
 \global\Widthtoks@=\expandafter{\the\Widthtoks@\fi\relax}%
 \edef\Width@##1##2{\i@=##1\relax\j@=##2\relax\the\Widthtoks@}%
 \global\Heighttoks@=\expandafter{\the\Heighttoks@\fi\relax}%
 \edef\Height@##1##2{\i@=##1\relax\j@=##2\relax\the\Heighttoks@}%
 \global\Depthtoks@=\expandafter{\the\Depthtoks@\fi\relax}%
 \edef\Depth@##1##2{\i@=##1\relax\j@=##2\relax\the\Depthtoks@}%
 \edef\next@{\the\Rowheighttoks@\noexpand\fi\relax}%
 \global\Rowheighttoks@=\expandafter{\next@}%
 \edef\Rowheight@##1{\i@=##1\relax\the\Rowheighttoks@}%
 \edef\next@{\the\Rowdepthtoks@\noexpand\fi\relax}%
 \global\Rowdepthtoks@=\expandafter{\next@}%
 \edef\Rowdepth@##1{\i@=##1\relax\the\Rowdepthtoks@}%
 \colwidthtoks@{\fi}%
 \setbox\zer@\vbox{%
  \unvbox\zer@
  \count@\rowcount@
  \loop
   \unskip\unpenalty
   \setbox\zer@\lastbox
   \ifnum\count@>\maxcolrow@ \advance\count@\m@ne
   \repeat
  \hbox{%
   \unhbox\zer@
   \count@\z@
   \loop
    \unskip
    \setbox\zer@\lastbox
    \edef\next@{\noexpand\or\noexpand\getdim@\the\wd\zer@\the\colwidthtoks@}%
     \global\colwidthtoks@=\expandafter{\next@}%
    \advance\count@\@ne
    \ifnum\count@<\Colcount@
    \repeat}}%
 \edef\next@{\noexpand\ifcase\noexpand\i@\the\colwidthtoks@}%
  \global\colwidthtoks@=\expandafter{\next@}%
 \edef\Colwidth@##1{\i@=##1\relax\the\colwidthtoks@}%
 \colwidthtoks@{}\Rowheighttoks@{}\Rowdepthtoks@{}\widthtoks@{}%
 \Widthtoks@{}\heighttoks@{}\Heighttoks@{}\depthtoks@{}\Depthtoks@{}%
}
\newcount\xoff@
\newcount\yoff@
\newcount\endcount@
\newcount\rcount@
\newdimen\firstx@
\newdimen\firsty@
\newdimen\secondx@
\newdimen\secondy@
\newdimen\tocenter@
\newdimen\charht@
\newdimen\charwd@
\def\outside@{\Err@{This arrow points outside the \string\newCD}}
\newif\ifsvertex@
\newif\iftvertex@
\def\arrow@#1#2{\xoff@=#1\relax\yoff@=#2\relax
 \count@\rowcount@ \advance\count@-\yoff@
 \ifnum\count@<\@ne \outside@ \else \ifnum\count@>\Rowcount@ \outside@ \fi\fi
 \count@\colcount@ \advance\count@\xoff@
 \ifnum\count@<\@ne \outside@ \else \ifnum\count@>\Colcount@ \outside@\fi\fi
 \tcolcount@\colcount@ \advance\tcolcount@\xoff@
 \Width@\rowcount@\colcount@ \tocenter@=-\getdim@ \divide\tocenter@\tw@
 \ifdim\getdim@=\z@
  \firstx@\z@ \firsty@\mathaxis@ \svertex@true
 \else
  \svertex@false
  \ifHshort@
   \Colwidth@\colcount@
    \ifE@ \firstx@=.5\getdim@ \else \firstx@=-.5\getdim@ \fi
  \else
   \ifE@ \firstx@=\getdim@ \else \firstx@=-\getdim@ \fi
   \divide\firstx@\tw@
  \fi
  \ifE@
   \ifH@ \advance\firstx@\thr@@\p@ \else \advance\firstx@-\thr@@\p@ \fi
  \else
   \ifH@ \advance\firstx@-\thr@@\p@ \else \advance\firstx@\thr@@\p@ \fi
  \fi
  \ifN@
   \Height@\rowcount@\colcount@ \firsty@=\getdim@
   \ifV@ \advance\firsty@\thr@@\p@ \fi
  \else
   \ifV@
    \Depth@\rowcount@\colcount@ \firsty@=-\getdim@
    \advance\firsty@-\thr@@\p@
   \else
    \firsty@\z@
   \fi
  \fi
 \fi
 \ifV@
 \else
  \Colwidth@\colcount@
  \ifE@ \secondx@=\getdim@ \else \secondx@=-\getdim@ \fi
  \divide\secondx@\tw@
  \ifE@ \else \getcgap@\colcount@ \advance\secondx@-\getdim@ \fi
  \endcount@=\colcount@ \advance\endcount@\xoff@
  \count@=\colcount@
  \ifE@
   \advance\count@\@ne
   \loop
    \ifnum\count@<\endcount@
    \Colwidth@\count@ \advance\secondx@\getdim@
    \getcgap@\count@ \advance\secondx@\getdim@
    \advance\count@\@ne
    \repeat
  \else
   \advance\count@\m@ne
   \loop
    \ifnum\count@>\endcount@
    \Colwidth@\count@ \advance\secondx@-\getdim@
    \getcgap@\count@ \advance\secondx@-\getdim@
    \advance\count@\m@ne
    \repeat
  \fi
  \Colwidth@\count@ \divide\getdim@\tw@
  \ifHshort@
  \else
   \ifE@ \advance\secondx@\getdim@ \else \advance\secondx@-\getdim@ \fi
  \fi
  \ifE@ \getcgap@\count@ \advance\secondx@\getdim@ \fi
  \rcount@\rowcount@ \advance\rcount@-\yoff@
  \Width@\rcount@\count@ \divide\getdim@\tw@
  \tvertex@false
  \ifH@\ifdim\getdim@=\z@\tvertex@true\Hshort@false\fi\fi
  \ifHshort@
  \else
   \ifE@ \advance\secondx@-\getdim@ \else \advance\secondx@\getdim@ \fi
  \fi
  \iftvertex@
   \advance\secondx@.4\p@
  \else
   \ifE@ \advance\secondx@-\thr@@\p@ \else \advance\secondx@\thr@@\p@ \fi
  \fi
 \fi
 \ifH@
 \else
  \ifN@
   \Rowheight@\rowcount@ \secondy@\getdim@
  \else
   \Rowdepth@\rowcount@ \secondy@-\getdim@
   \getrgap@\rowcount@ \advance\secondy@-\getdim@
  \fi
  \endcount@=\rowcount@ \advance\endcount@-\yoff@
  \count@=\rowcount@
  \ifN@
   \advance\count@\m@ne
   \loop
    \ifnum\count@>\endcount@
    \Rowheight@\count@ \advance\secondy@\getdim@
    \Rowdepth@\count@ \advance\secondy@\getdim@
    \getrgap@\count@ \advance\secondy@\getdim@
    \advance\count@\m@ne
    \repeat
  \else
   \advance\count@\@ne
   \loop
    \ifnum\count@<\endcount@
    \Rowheight@\count@ \advance\secondy@-\getdim@
    \Rowdepth@\count@ \advance\secondy@-\getdim@
    \getrgap@\count@ \advance\secondy@-\getdim@
    \advance\count@\@ne
    \repeat
  \fi
  \tvertex@false
  \ifV@\Width@\count@\colcount@\ifdim\getdim@=\z@\tvertex@true\fi\fi
  \ifN@
   \getrgap@\count@ \advance\secondy@\getdim@
   \Rowdepth@\count@ \advance\secondy@\getdim@
   \iftvertex@
    \advance\secondy@\mathaxis@
   \else
    \Depth@\count@\tcolcount@ \advance\secondy@-\getdim@
    \advance\secondy@-\thr@@\p@
   \fi
  \else
   \Rowheight@\count@ \advance\secondy@-\getdim@
   \iftvertex@
    \advance\secondy@\mathaxis@
   \else
    \Height@\count@\tcolcount@ \advance\secondy@\getdim@
    \advance\secondy@\thr@@\p@
   \fi
  \fi
 \fi
 \ifV@\else\advance\firstx@\sxdimen@\fi
 \ifH@\else\advance\firsty@\sydimen@\fi
 \iftX@
  \advance\secondy@\tXdimen@ii
  \advance\secondx@\tXdimen@i
  \slope@
 \else
  \iftY@
   \advance\secondy@\tYdimen@ii
   \advance\secondx@\tYdimen@i
   \slope@
   \secondy@=\secondx@ \advance\secondy@-\firstx@
   \ifNESW@ \else \multiply\secondy@\m@ne \fi
   \multiply\secondy@\tan@i \divide\secondy@\tan@ii \advance\secondy@\firsty@
  \else
   \ifa@
    \slope@
    \ifNESW@ \global\advance\angcount@\exacount@ \else
      \global\advance\angcount@-\exacount@ \fi
    \ifnum\angcount@>23 \angcount@23 \fi
    \ifnum\angcount@<\@ne \angcount@\@ne \fi
    \slope@a\angcount@
    \ifY@
     \advance\secondy@\Ydimen@
    \else
     \ifX@
      \advance\secondx@\Xdimen@
      \dimen@\secondx@ \advance\dimen@-\firstx@
      \ifNESW@\else\multiply\dimen@\m@ne\fi
      \multiply\dimen@\tan@i \divide\dimen@\tan@ii
      \advance\dimen@\firsty@ \secondy@=\dimen@
     \fi
    \fi
   \else
    \ifH@\else\ifV@\else\slope@\fi\fi
   \fi
  \fi
 \fi
 \ifH@\else\ifV@\else\ifsvertex@\else
  \dimen@=6\p@ \multiply\dimen@\tan@ii
  \count@=\tan@i \advance\count@\tan@ii \divide\dimen@\count@
  \ifE@ \advance\firstx@\dimen@ \else \advance\firstx@-\dimen@ \fi
  \multiply\dimen@\tan@i \divide\dimen@\tan@ii
  \ifN@ \advance\firsty@\dimen@ \else \advance\firsty@-\dimen@ \fi
 \fi\fi\fi
 \ifp@
  \ifH@\else\ifV@\else
   \getcos@\pdimen@ \advance\firsty@\dimen@ \advance\secondy@\dimen@
   \ifNESW@ \advance\firstx@-\dimen@ii \else \advance\firstx@\dimen@ii \fi
  \fi\fi
 \fi
 \ifH@\else\ifV@\else
  \ifnum\tan@i>\tan@ii
   \charht@=10\p@ \charwd@=10\p@
   \multiply\charwd@\tan@ii \divide\charwd@\tan@i
  \else
   \charwd@=10\p@ \charht@=10\p@
   \divide\charht@\tan@ii \multiply\charht@\tan@i
  \fi
  \ifnum\tcount@=\thr@@
   \ifN@ \advance\secondy@-.3\charht@ \else\advance\secondy@.3\charht@ \fi
  \fi
  \ifnum\scount@=\tw@
   \ifE@ \advance\firstx@.3\charht@ \else \advance\firstx@-.3\charht@ \fi
  \fi
  \ifnum\tcount@=12
   \ifN@ \advance\secondy@-\charht@ \else \advance\secondy@\charht@ \fi
  \fi
  \iftY@
  \else
   \ifa@
    \ifX@
    \else
     \secondx@\secondy@ \advance\secondx@-\firsty@
     \ifNESW@\else\multiply\secondx@\m@ne\fi
     \multiply\secondx@\tan@ii \divide\secondx@\tan@i
     \advance\secondx@\firstx@
    \fi
   \fi
  \fi
 \fi\fi
 \ifH@\harrow@\else\ifV@\varrow@\else\arrow@@\fi\fi}
\newdimen\mathaxis@
\mathaxis@90\p@ \divide\mathaxis@36
\def\harrow@b{\ifE@\hskip\tocenter@\hskip\firstx@\fi}
\def\harrow@bb{\ifE@\hskip\xdimen@\else\hskip\Xdimen@\fi}
\def\harrow@e{\ifE@\else\hskip-\firstx@\hskip-\tocenter@\fi}
\def\harrow@ee{\ifE@\hskip-\Xdimen@\else\hskip-\xdimen@\fi}
\def\harrow@{\dimen@\secondx@\advance\dimen@-\firstx@
 \ifE@ \let\next@\rlap \else  \multiply\dimen@\m@ne \let\next@\llap \fi
 \next@{%
  \harrow@b
  \smash{\raise\pdimen@\hbox to\dimen@
   {\harrow@bb\arrow@ii
    \ifnum\arrcount@=\m@ne \else \ifnum\arrcount@=\thr@@ \else
     \ifE@
      \ifnum\scount@=\m@ne
      \else
       \ifcase\scount@\or\or\char118 \or\char117 \or\or\or\char119 \or
       \char120 \or\char121 \or\char122 \or\or\or\arrow@i\char125 \or
       \char117 \hskip\thr@@\p@\char117 \hskip-\thr@@\p@\fi
      \fi
     \else
      \ifnum\tcount@=\m@ne
      \else
       \ifcase\tcount@\char117 \or\or\char117 \or\char118 \or\char119 \or
       \char120\or\or\or\or\or\char121 \or\char122 \or\arrow@i\char125
       \or\char117 \hskip\thr@@\p@\char117 \hskip-\thr@@\p@\fi
      \fi
     \fi
    \fi\fi
    \dimen@\mathaxis@ \advance\dimen@.2\p@
    \dimen@ii\mathaxis@ \advance\dimen@ii-.2\p@
    \ifnum\arrcount@=\m@ne
     \let\leads@\null
    \else
     \ifcase\arrcount@
      \def\leads@{\hrule height\dimen@ depth-\dimen@ii}\or
      \def\leads@{\hrule height\dimen@ depth-\dimen@ii}\or
      \def\leads@{\hbox to10\p@{%
       \leaders\hrule height\dimen@ depth-\dimen@ii\hfil
       \hfil
      \leaders\hrule height\dimen@ depth-\dimen@ii\hskip\z@ plus2fil\relax
       \hfil
       \leaders\hrule height\dimen@ depth-\dimen@ii\hfil}}\or
     \def\leads@{\hbox{\hbox to10\p@{\dimen@\mathaxis@ \advance\dimen@1.2\p@
       \dimen@ii\dimen@ \advance\dimen@ii-.4\p@
       \leaders\hrule height\dimen@ depth-\dimen@ii\hfil}%
       \kern-10\p@
       \hbox to10\p@{\dimen@\mathaxis@ \advance\dimen@-1.2\p@
       \dimen@ii\dimen@ \advance\dimen@ii-.4\p@
       \leaders\hrule height\dimen@ depth-\dimen@ii\hfil}}}\fi
    \fi
    \cleaders\leads@\hfil
    \ifnum\arrcount@=\m@ne\else\ifnum\arrcount@=\thr@@\else
     \arrow@i
     \ifE@
      \ifnum\tcount@=\m@ne
      \else
       \ifcase\tcount@\char119 \or\or\char119 \or\char120 \or\char121 \or
       \char122 \or \or\or\or\or\char123\or\char124 \or
       \char125 \or\char119 \hskip-\thr@@\p@\char119 \hskip\thr@@\p@\fi
      \fi
     \else
      \ifcase\scount@\or\or\char120 \or\char119 \or\or\or\char121 \or\char122
      \or\char123 \or\char124 \or\or\or\char125 \or
      \char119 \hskip-\thr@@\p@\char119 \hskip\thr@@\p@\fi
     \fi
    \fi\fi
    \harrow@ee}}%
  \harrow@e}%
 \iflabel@i
  \dimen@ii\z@ \setbox\zer@\hbox{$\m@th\tsize@@\label@i$}%
  \ifnum\arrcount@=\m@ne
  \else
   \advance\dimen@ii\mathaxis@
   \advance\dimen@ii\dp\zer@ \advance\dimen@ii\tw@\p@
   \ifnum\arrcount@=\thr@@ \advance\dimen@ii\tw@\p@ \fi
  \fi
  \advance\dimen@ii\pdimen@
  \next@{\harrow@b\smash{\raise\dimen@ii\hbox to\dimen@
   {\harrow@bb\hskip\tw@\ldimen@i\hfil\box\zer@\hfil\harrow@ee}}\harrow@e}%
 \fi
 \iflabel@ii
  \ifnum\arrcount@=\m@ne
  \else
   \setbox\zer@\hbox{$\m@th\tsize@\label@ii$}%
   \dimen@ii-\ht\zer@ \advance\dimen@ii-\tw@\p@
   \ifnum\arrcount@=\thr@@ \advance\dimen@ii-\tw@\p@ \fi
   \advance\dimen@ii\mathaxis@ \advance\dimen@ii\pdimen@
   \next@{\harrow@b\smash{\raise\dimen@ii\hbox to\dimen@
    {\harrow@bb\hskip\tw@\ldimen@ii\hfil\box\zer@\hfil\harrow@ee}}\harrow@e}%
  \fi
 \fi}
\let\tsize@\tsize
\def\tsizenewCDlabels{\let\tsize@\tsize}
\def\ssizenewCDlabels{\let\tsize@\ssize}
\def\tsize@@{\ifnum\arrcount@=\m@ne\else\tsize@\fi}
\def\varrow@{\dimen@\secondy@ \advance\dimen@-\firsty@
 \ifN@ \else \multiply\dimen@\m@ne \fi
 \setbox\zer@\vbox to\dimen@
  {\ifN@ \vskip-\Ydimen@ \else \vskip\ydimen@ \fi
   \ifnum\arrcount@=\m@ne\else\ifnum\arrcount@=\thr@@\else
    \hbox{\arrow@iii
     \ifN@
      \ifnum\tcount@=\m@ne
      \else
       \ifcase\tcount@\char117 \or\or\char117 \or\char118 \or\char119 \or
       \char120 \or\or\or\or\or\char121 \or\char122 \or\char123 \or
       \vbox{\hbox{\char117 }\nointerlineskip\vskip\thr@@\p@
       \hbox{\char117 }\vskip-\thr@@\p@}\fi
      \fi
     \else
      \ifcase\scount@\or\or\char118 \or\char117 \or\or\or\char119 \or
      \char120 \or\char121 \or\char122 \or\or\or\char123 \or
      \vbox{\hbox{\char117 }\nointerlineskip\vskip\thr@@\p@
      \hbox{\char117 }\vskip-\thr@@\p@}\fi
     \fi}%
    \nointerlineskip
   \fi\fi
   \ifnum\arrcount@=\m@ne
    \let\leads@\null
   \else
    \ifcase\arrcount@\let\leads@\vrule\or\let\leads@\vrule\or
    \def\leads@{\vbox to10\p@{%
     \hrule height 1.67\p@ depth\z@ width.4\p@
     \vfil
     \hrule height 3.33\p@ depth\z@ width.4\p@
     \vfil
     \hrule height 1.67\p@ depth\z@ width.4\p@}}\or
    \def\leads@{\hbox{\vrule height\p@\hskip\tw@\p@\vrule}}\fi
   \fi
  \cleaders\leads@\vfill\nointerlineskip
   \ifnum\arrcount@=\m@ne\else\ifnum\arrcount@=\thr@@\else
    \hbox{\arrow@iv
     \ifN@
      \ifcase\scount@\or\or\char118 \or\char117 \or\or\or\char119 \or
      \char120 \or\char121 \or\char122 \or\or\or\arrow@iii\char123 \or
      \vbox{\hbox{\char117 }\nointerlineskip\vskip-\thr@@\p@
      \hbox{\char117 }\vskip\thr@@\p@}\fi
     \else
      \ifnum\tcount@=\m@ne
      \else
       \ifcase\tcount@\char117 \or\or\char117 \or\char118 \or\char119 \or
       \char120 \or\or\or\or\or\char121 \or\char122 \or\arrow@iii\char123 \or
       \vbox{\hbox{\char117 }\nointerlineskip\vskip-\thr@@\p@
       \hbox{\char117 }\vskip\thr@@\p@}\fi
      \fi
     \fi}%
   \fi\fi
   \ifN@\vskip\ydimen@\else\vskip-\Ydimen@\fi}%
 \ifN@
  \dimen@ii\firsty@
 \else
  \dimen@ii-\firsty@ \advance\dimen@ii\ht\zer@ \multiply\dimen@ii\m@ne
 \fi
 \rlap{\smash{\hskip\tocenter@ \hskip\pdimen@ \raise\dimen@ii \box\zer@}}%
 \iflabel@i
  \setbox\zer@\vbox to\dimen@{\vfil
   \hbox{$\m@th\tsize@@\label@i$}\vskip\tw@\ldimen@i\vfil}%
  \rlap{\smash{\hskip\tocenter@ \hskip\pdimen@
  \ifnum\arrcount@=\m@ne \let\next@\relax \else \let\next@\llap \fi
  \next@{\raise\dimen@ii\hbox{\ifnum\arrcount@=\m@ne \hskip-.5\wd\zer@ \fi
   \box\zer@ \ifnum\arrcount@=\m@ne \else \hskip\tw@\p@ \fi}}}}%
 \fi
 \iflabel@ii
  \ifnum\arrcount@=\m@ne
  \else
   \setbox\zer@\vbox to\dimen@{\vfil
    \hbox{$\m@th\tsize@\label@ii$}\vskip\tw@\ldimen@ii\vfil}%
   \rlap{\smash{\hskip\tocenter@ \hskip\pdimen@
   \rlap{\raise\dimen@ii\hbox{\ifnum\arrcount@=\thr@@ \hskip4.5\p@ \else
    \hskip2.5\p@ \fi\box\zer@}}}}%
  \fi
 \fi
}
\newdimen\goal@
\newdimen\shifted@
\newcount\Tcount@
\newcount\Scount@
\newbox\shaft@
\newcount\slcount@
\def\getcos@#1{%
 \ifnum\tan@i<\tan@ii
  \dimen@#1%
  \ifnum\slcount@<8 \count@9 \else \ifnum\slcount@<12 \count@8 \else
   \count@7 \fi\fi
  \multiply\dimen@\count@ \divide\dimen@10
  \dimen@ii\dimen@ \multiply\dimen@ii\tan@i \divide\dimen@ii\tan@ii
 \else
  \dimen@ii#1%
  \count@-\slcount@ \advance\count@24
  \ifnum\count@<8 \count@9 \else \ifnum\count@<12 \count@8
   \else\count@7 \fi\fi
  \multiply\dimen@ii\count@ \divide\dimen@ii10
  \dimen@\dimen@ii \multiply\dimen@\tan@ii \divide\dimen@\tan@i
 \fi}
\newdimen\adjust@
\def\Nnext@{\ifN@\let\next@\raise\else\let\next@\lower\fi}
\def\arrow@@{\slcount@\angcount@
 \ifNESW@
  \ifnum\angcount@<10
   \let\arrowfont@=\arrow@i \advance\angcount@\m@ne \multiply\angcount@13
  \else
   \ifnum\angcount@<19
    \let\arrowfont@=\arrow@ii \advance\angcount@-10 \multiply\angcount@13
   \else
    \let\arrowfont@=\arrow@iii \advance\angcount@-19 \multiply\angcount@13
  \fi\fi
  \Tcount@\angcount@
 \else
  \ifnum\angcount@<5
   \let\arrowfont@=\arrow@iii \advance\angcount@\m@ne \multiply\angcount@13
   \advance\angcount@65
  \else
   \ifnum\angcount@<14
    \let\arrowfont@=\arrow@iv \advance\angcount@-5 \multiply\angcount@13
   \else
    \ifnum\angcount@<23
     \let\arrowfont@=\arrow@v \advance\angcount@-14 \multiply\angcount@13
    \else
     \let\arrowfont@=\arrow@i \angcount@=117
  \fi\fi\fi
  \ifnum\angcount@=117 \Tcount@=115 \else\Tcount@\angcount@ \fi
 \fi
 \Scount@\Tcount@
 \ifE@
  \ifnum\tcount@=\z@ \advance\Tcount@\tw@ \else\ifnum\tcount@=13
   \advance\Tcount@\tw@ \else \advance\Tcount@\tcount@ \fi\fi
  \ifnum\scount@=\z@ \else \ifnum\scount@=13 \advance\Scount@\thr@@ \else
   \advance\Scount@\scount@ \fi\fi
 \else
  \ifcase\tcount@\advance\Tcount@\thr@@\or\or\advance\Tcount@\thr@@\or
  \advance\Tcount@\tw@\or\advance\Tcount@6 \or\advance\Tcount@7
  \or\or\or\or\or \advance\Tcount@8 \or\advance\Tcount@9 \or
  \advance\Tcount@12 \or\advance\Tcount@\thr@@\fi
  \ifcase\scount@\or\or\advance\Scount@\thr@@\or\advance\Scount@\tw@\or
  \or\or\advance\Scount@4 \or\advance\Scount@5 \or\advance\Scount@10
  \or\advance\Scount@11 \or\or\or\advance\Scount@12 \or\advance
  \Scount@\tw@\fi
 \fi
 \ifcase\arrcount@\or\or\advance\angcount@\@ne\else\fi
 \ifN@ \shifted@=\firsty@ \else\shifted@=-\firsty@ \fi
 \ifE@ \else\advance\shifted@\charht@ \fi
 \goal@=\secondy@ \advance\goal@-\firsty@
 \ifN@\else\multiply\goal@\m@ne\fi
 \setbox\shaft@\hbox{\arrowfont@\char\angcount@}%
 \ifnum\arrcount@=\thr@@
  \getcos@{1.5\p@}%
  \setbox\shaft@\hbox to\wd\shaft@{\arrowfont@
   \rlap{\hskip\dimen@ii
    \smash{\ifNESW@\let\next@\lower\else\let\next@\raise\fi
     \next@\dimen@\hbox{\arrowfont@\char\angcount@}}}%
   \rlap{\hskip-\dimen@ii
    \smash{\ifNESW@\let\next@\raise\else\let\next@\lower\fi
      \next@\dimen@\hbox{\arrowfont@\char\angcount@}}}\hfil}%
 \fi
 \rlap{\smash{\hskip\tocenter@\hskip\firstx@
  \ifnum\arrcount@=\m@ne
  \else
   \ifnum\arrcount@=\thr@@
   \else
    \ifnum\scount@=\m@ne
    \else
     \ifnum\scount@=\z@
     \else
      \setbox\zer@\hbox{\ifnum\angcount@=117 \arrow@v\else\arrowfont@\fi
       \char\Scount@}%
      \ifNESW@
       \ifnum\scount@=\tw@
        \dimen@=\shifted@ \advance\dimen@-\charht@
        \ifN@\hskip-\wd\zer@\fi
        \Nnext@
        \next@\dimen@\copy\zer@
        \ifN@\else\hskip-\wd\zer@\fi
       \else
        \Nnext@
        \ifN@\else\hskip-\wd\zer@\fi
        \next@\shifted@\copy\zer@
        \ifN@\hskip-\wd\zer@\fi
       \fi
       \ifnum\scount@=12
        \advance\shifted@\charht@ \advance\goal@-\charht@
        \ifN@ \hskip\wd\zer@ \else \hskip-\wd\zer@ \fi
       \fi
       \ifnum\scount@=13
        \getcos@{\thr@@\p@}%
        \ifN@ \hskip\dimen@ \else \hskip-\wd\zer@ \hskip-\dimen@ \fi
        \adjust@\shifted@ \advance\adjust@\dimen@ii
        \Nnext@
        \next@\adjust@\copy\zer@
        \ifN@ \hskip-\dimen@ \hskip-\wd\zer@ \else \hskip\dimen@ \fi
       \fi
      \else
       \ifN@\hskip-\wd\zer@\fi
       \ifnum\scount@=\tw@
        \ifN@ \hskip\wd\zer@ \else \hskip-\wd\zer@ \fi
        \dimen@=\shifted@ \advance\dimen@-\charht@
        \Nnext@
        \next@\dimen@\copy\zer@
        \ifN@\hskip-\wd\zer@\fi
       \else
        \Nnext@
        \next@\shifted@\copy\zer@
        \ifN@\else\hskip-\wd\zer@\fi
       \fi
       \ifnum\scount@=12
        \advance\shifted@\charht@ \advance\goal@-\charht@
        \ifN@ \hskip-\wd\zer@ \else \hskip\wd\zer@ \fi
       \fi
       \ifnum\scount@=13
        \getcos@{\thr@@\p@}%
        \ifN@ \hskip-\wd\zer@ \hskip-\dimen@ \else \hskip\dimen@ \fi
        \adjust@\shifted@ \advance\adjust@\dimen@ii
        \Nnext@
        \next@\adjust@\copy\zer@
        \ifN@ \hskip\dimen@ \else \hskip-\dimen@ \hskip-\wd\zer@ \fi
       \fi	
      \fi
  \fi\fi\fi\fi
  \ifnum\arrcount@=\m@ne
  \else
   \loop
    \ifdim\goal@>\charht@
    \ifE@\else\hskip-\charwd@\fi
    \Nnext@
    \next@\shifted@\copy\shaft@
    \ifE@\else\hskip-\charwd@\fi
    \advance\shifted@\charht@ \advance\goal@ -\charht@
    \repeat
   \ifdim\goal@>\z@
    \dimen@=\charht@ \advance\dimen@-\goal@
    \divide\dimen@\tan@i \multiply\dimen@\tan@ii
    \ifE@ \hskip-\dimen@ \else \hskip-\charwd@ \hskip\dimen@ \fi
    \adjust@=\shifted@ \advance\adjust@-\charht@ \advance\adjust@\goal@
    \Nnext@
    \next@\adjust@\copy\shaft@
    \ifE@ \else \hskip-\charwd@ \fi
   \else
    \adjust@=\shifted@ \advance\adjust@-\charht@
   \fi
  \fi
  \ifnum\arrcount@=\m@ne
  \else
   \ifnum\arrcount@=\thr@@
   \else
    \ifnum\tcount@=\m@ne
    \else
     \setbox\zer@
      \hbox{\ifnum\angcount@=117 \arrow@v\else\arrowfont@\fi\char\Tcount@}%
     \ifnum\tcount@=\thr@@
      \advance\adjust@\charht@
      \ifE@\else\ifN@\hskip-\charwd@\else\hskip-\wd\zer@\fi\fi
     \else
      \ifnum\tcount@=12
       \advance\adjust@\charht@
       \ifE@\else\ifN@\hskip-\charwd@\else\hskip-\wd\zer@\fi\fi
      \else
       \ifE@\hskip-\wd\zer@\fi
     \fi\fi
     \Nnext@
     \next@\adjust@\copy\zer@
     \ifnum\tcount@=13
      \hskip-\wd\zer@
      \getcos@{\thr@@\p@}%
      \ifE@\hskip-\dimen@ \else\hskip\dimen@ \fi
      \advance\adjust@-\dimen@ii
      \Nnext@
      \next@\adjust@\box\zer@
     \fi
  \fi\fi\fi}}%
 \iflabel@i
  \rlap{\hskip\tocenter@
  \dimen@\firstx@ \advance\dimen@\secondx@ \divide\dimen@\tw@
  \advance\dimen@\ldimen@i
  \dimen@ii\firsty@ \advance\dimen@ii\secondy@ \divide\dimen@ii\tw@
  \multiply\ldimen@i\tan@i \divide\ldimen@i\tan@ii
  \ifNESW@ \advance\dimen@ii\ldimen@i \else \advance\dimen@ii-\ldimen@i \fi
  \setbox\zer@\hbox{\ifNESW@\else\ifnum\arrcount@=\thr@@\hskip4\p@\else
   \hskip\tw@\p@\fi\fi
   $\m@th\tsize@@\label@i$\ifNESW@\ifnum\arrcount@=\thr@@\hskip4\p@\else
   \hskip\tw@\p@\fi\fi}%
  \ifnum\arrcount@=\m@ne
   \ifNESW@ \advance\dimen@.5\wd\zer@ \advance\dimen@\p@ \else
    \advance\dimen@-.5\wd\zer@ \advance\dimen@-\p@ \fi
   \advance\dimen@ii-.5\ht\zer@
  \else
   \advance\dimen@ii\dp\zer@
   \ifnum\slcount@<6 \advance\dimen@ii\tw@\p@ \fi
  \fi
  \hskip\dimen@
  \ifNESW@ \let\next@\llap \else\let\next@\rlap \fi
  \next@{\smash{\raise\dimen@ii\box\zer@}}}%
 \fi
 \iflabel@ii
  \ifnum\arrcount@=\m@ne
  \else
   \rlap{\hskip\tocenter@
   \dimen@\firstx@ \advance\dimen@\secondx@ \divide\dimen@\tw@
   \ifNESW@ \advance\dimen@\ldimen@ii \else \advance\dimen@-\ldimen@ii \fi
   \dimen@ii\firsty@ \advance\dimen@ii\secondy@ \divide\dimen@ii\tw@
   \multiply\ldimen@ii\tan@i \divide\ldimen@ii\tan@ii
   \advance\dimen@ii\ldimen@ii
   \setbox\zer@\hbox{\ifNESW@\ifnum\arrcount@=\thr@@\hskip4\p@\else
    \hskip\tw@\p@\fi\fi
    $\m@th\tsize@\label@ii$\ifNESW@\else\ifnum\arrcount@=\thr@@\hskip4\p@
    \else\hskip\tw@\p@\fi\fi}%
   \advance\dimen@ii-\ht\zer@
   \ifnum\slcount@<9 \advance\dimen@ii-\thr@@\p@ \fi
   \ifNESW@ \let\next@\rlap \else \let\next@\llap \fi
   \hskip\dimen@\next@{\smash{\raise\dimen@ii\box\zer@}}}%
  \fi
 \fi
}
\def\outnewCD@#1{\def#1{\Err@{\string#1 must not be used within \string\newCD}}}
\newskip\prenewCDskip@
\newskip\postnewCDskip@
\prenewCDskip@\z@
\postnewCDskip@\z@
\def\prenewCDspace#1{\RIfMIfI@
 \onlydmatherr@\prenewCDspace\else\advance\prenewCDskip@#1\relax\fi\else
 \onlydmatherr@\prenewCDspace\fi}
\def\postnewCDspace#1{\RIfMIfI@
 \onlydmatherr@\postnewCDspace\else\advance\postnewCDskip@#1\relax\fi\else
 \onlydmatherr@\postnewCDspace\fi}
\def\predisplayspace#1{\RIfMIfI@
 \onlydmatherr@\predisplayspace\else
 \advance\abovedisplayskip#1\relax
 \advance\abovedisplayshortskip#1\relax\fi
 \else\onlydmatherr@\prenewCDspace\fi}
\def\postdisplayspace#1{\RIfMIfI@
 \onlydmatherr@\postdisplayspace\else
 \advance\belowdisplayskip#1\relax
 \advance\belowdisplayshortskip#1\relax\fi
 \else\onlydmatherr@\postdisplayspace\fi}
\def\PrenewCDSpace#1{\global\prenewCDskip@#1\relax}
\def\PostnewCDSpace#1{\global\postnewCDskip@#1\relax}
\def\newCD#1\endnewCD{%
 \outnewCD@\cgaps\outnewCD@\rgaps\outnewCD@\Cgaps\outnewCD@\Rgaps
 \prenewCD@#1\endnewCD
 \advance\abovedisplayskip\prenewCDskip@
 \advance\abovedisplayshortskip\prenewCDskip@
 \advance\belowdisplayskip\postnewCDskip@
 \advance\belowdisplayshortskip\postnewCDskip@
 \vcenter{\vskip\prenewCDskip@ \Let@ \colcount@\@ne \rowcount@\z@
  \everycr{%
   \noalign{%
    \ifnum\rowcount@=\Rowcount@
    \else
     \global\nointerlineskip
     \getrgap@\rowcount@ \vskip\getdim@
     \global\advance\rowcount@\@ne \global\colcount@\@ne
    \fi}}%
  \tabskip\z@
  \halign{&\global\xoff@\z@ \global\yoff@\z@
   \getcgap@\colcount@ \hskip\getdim@
   \hfil\vrule height10\p@ width\z@ depth\z@
   $\m@th\displaystyle{##}$\hfil
   \global\advance\colcount@\@ne\cr
   #1\crcr}\vskip\postnewCDskip@}%
 \prenewCDskip@\z@\postnewCDskip@\z@
 \def\getcgap@##1{\ifcase##1\or\getdim@\z@\else\getdim@\standardcgap\fi}%
 \def\getrgap@##1{\ifcase##1\getdim@\z@\else\getdim@\standardrgap\fi}%
 \let\Width@\relax\let\Height@\relax\let\Depth@\relax\let\Rowheight@\relax
 \let\Rowdepth@\relax\let\Colwdith@\relax
}
\catcode`\@=\active
\input amsppt.sty
\hsize 30pc
\vsize 47pc
\magnification=\magstep1
\def\nmb#1#2{#2}         
\def\cit#1#2{\ifx#1!\cite{#2}\else#2\fi} 
\def\totoc{}             
\def\idx{}               
\def\ign#1{}             
\redefine\o{\circ}
\define\X{\frak X}
\define\al{\alpha}
\define\be{\beta}

\define\de{\delta}
\define\ep{\varepsilon}

\define\th{\theta}

\define\la{\lambda}

\define\ta{\tau}
\define\ph{\varphi}

\define\ps{\psi}
\define\om{\omega}
\define\Ga{\Gamma}

\define\Ph{\Phi}
\define\Ps{\Psi}
\define\Om{\Omega}
\predefine\ii\i
\redefine\i{^{-1}}
\define\x{\times}

\define\Der{{\operatorname{Der}}}

\define\ad{\operatorname{ad}}

\define\Lip{\operatorname{Lip}}
\define\Int{\operatorname{Inn}}
\define\Out{{\operatorname{Out}}}

\define\Id{{\operatorname{Id}}}
\let\on\operatorname
\redefine\L{{\Cal L}}
\def\today{\ifcase\month\or
 January\or February\or March\or April\or May\or June\or
 July\or August\or September\or October\or November\or December\fi
 \space\number\day, \number\year}
\hyphenation{ho-mo-mor-phism boun-ded}

\topmatter
\title  
Smooth $*$-Algebras 
\endtitle
\author  
Michel Dubois-Violette\\ 
Andreas Kriegl\\
Yoshiaki Maeda\\
Peter W. Michor  
\endauthor 
\leftheadtext{\smc M\. Dubois-Violette, A. Kriegl, Y\. Maeda, P\. Michor} 
\rightheadtext{\smc Smooth $*$-algebras} 
\address 
M\. Dubois-Violette:  
Laboratoire de Physique Th\'eorique et Hautes Energies,  
Universit\'e Paris XI, B\^atiment 211, F-91405 Orsay Cedex,  
France 
\endaddress 
\email flad\@qcd.th.u-psud.fr \endemail 
\address
A\. Kriegl: Institut f\"ur Mathematik, Universit\"at Wien,
Strudlhofgasse 4, A-1090 Wien, Austria
\endaddress
\address Y. Maeda: Department of Mathematics, 
Keio University, 3-14-1 Hiyoshi, Kohoku-ku, Yokohama 2238522, Japan
\endaddress
\email maeda\@math.keio.ac.jp \endemail
\address 
P\. Michor: Institut f\"ur Mathematik, Universit\"at Wien, 
Strudlhofgasse 4, A-1090 Wien, Austria; and:  
Erwin Schr\"odinger International Institute of Mathematical Physics,  
Boltzmanngasse 9, A-1090 Wien, Austria 
\endaddress 
\email Peter.Michor\@esi.ac.at \endemail 
\date {\today} \enddate 
\thanks  
P. Michor was supported  
by `Fonds zur F\"orderung der wissenschaftlichen  
Forschung, Projekt P~14195~MAT'.
\endthanks 
\abstract Looking for the universal covering of the smooth 
non-commutative torus leads to a curve of associative multiplications 
on the space $\Cal O_M'(\Bbb R^{2n})\cong \Cal O_C(\Bbb R^{2n})$ of 
Laurent Schwartz which is smooth in the deformation parameter $\hbar$. 
The Taylor expansion in $\hbar$ leads to the formal Moyal star 
product. The non-commutative torus and this version of the Heisenberg 
plane are examples of smooth *-algebras: smooth in the sense of 
having many derivations. A tentative definition of this concept is given. 
\endabstract

\endtopmatter

\document

\heading Table of contents \endheading
0. Introduction 

1. Smooth $*$-algebras 

2. The non-commutative torus 

3. The smooth Heisenberg algebra 

4. Appendix: Calculus in infinite dimensions and convenient 
vector spaces 

\head\totoc\nmb0{0}. Introduction \endhead

The noncommutative torus in its topological version 
($C^\ast$-completion) as well as in its smooth version \cit!{6} is one 
of the most important examples in noncommutative geometry. Beside the 
fact that the classical tools of differential geometry have 
unambiguous generalizations to it, it provides a very nontrivial 
example of noncommutative geometry satisfying the axioms of \cit!{7} 
(see also in \cit!{8}, \cit!{9}). We looked at its smooth version and 
asked for its universal covering. We found the Heisenberg plane as it 
is presented in this paper: a twisted convolution on a carefully 
chosen space of distributions, namely the topological dual space 
${\Cal O}'_{M}$ of the Schwartz space ${\Cal O}_{M}$ of smooth slowly 
increasing functions at $\infty$, \cit!{29}, \cit!{30}. 
It is large enough 
to contain the space of rapidly decreasing measures with support in 
the lattice $(2\pi\Bbb Z)^2$ that is a space isomorphic to the space 
of smooth functions on the noncommutative torus (as well as on the 
usual commutative torus). The multiplication turns out to be a smooth 
curve in the deformation parameter $\hbar$. Moreover, looking at it 
via Fourier transform, Taylor expansion of the multiplication in the 
deformation parameter $\hbar$ leads to the formal Moyal star-product 
which is well known from deformation quantization, 
\cit!{24}, \cit!{1}.
 
Then we noticed that we found examples of noncommutative 
$\ast$-algebras generalizing algebras of complex smooth  functions. 
These $\ast$-algebras which can be realized as $\ast$-algebras of 
unbounded operators in Hilbert space admit ``many" derivations 
specifying thereby the generalized smooth structure (see below). These 
algebras are defined in Section \nmb!{1} and are tentatively called smooth 
$\ast$-algebras.
 
Section \nmb!{2} contains our treatment of the smooth non-commutative 
torus, and also some related material like the smooth non-commutative 
circle of rational slope $b/a$, a quotient of the smooth 
non-commutative torus. 

The appendix in section \nmb!{4} gives an overview on convenient 
calculus in infinite dimensions which is necessary to obtain our 
results about smoothness in the deformation parameter $\hbar$, and 
which also gives the right setting for multilinear algebra with 
locally convex vector spaces.  

Work on this paper started in 1996, but we were unable to prove that 
the Heisenberg plane is a smooth *-algebra. Finally we gave up and 
stated this as a conjecture. The problem is finding enough states. 

\head\totoc\nmb0{1}. Smooth $*$-algebras \endhead

\subhead\nmb.{1.1}. Preliminaries \endsubhead
Throughout this paper by a $\ast$-{\it algebra} we always mean a
complex associative algebra $A$ with unit equipped with an antilinear
involution $f\mapsto f^\ast$ which reverses the order of products i.e.
which satisfies $(fg)^\ast=g^\ast f^\ast$, $\forall f,g\in A$. Given
a $\ast$-algebra $A$, a {\it hermitian representation} \cit!{25}
of $A$ in a Hilbert space ${\Cal H}$ is a homomorphism $\pi$ of unital
algebras of $A$ into the algebra of endomorphisms of a dense subspace
$D(\pi)$ of ${\Cal H}$ satisfying
$(\Ps,\pi(f)\Ph)=(\pi(f^\ast)\Ps,\Ph)$ for any $f\in A$ and
$\Ps,\Ph\in D(\pi)$; the dense subspace $D(\pi)$ of ${\Cal H}$
is refered to as {\it the domain of} $\pi$. The image
of a hermitian representation in $\Cal H$ is a unital subalgebra of
the algebra of endomorphisms of the dense domain $D$ of the
representation which is also a $\ast$-algebra for an obvious involution;
such a $\ast$-algebra will be refered to as a $\ast$-{\it algebra of
(unbounded) operators in the Hilbert space $\Cal H$ with domain $D$}.

A linear form $\ph$ on a $\ast$-algebra $A$ is said to be {\it 
positive} if $\ph(f^\ast f)\geq 0$ for all $f\in A$. Such a positive
linear form satisfies $\ph(f^\ast)=\overline{\ph(f)}$ (for all $f\in
A$) and $(f,g)_{\om}=\ph(f^\ast g)$ is a pre-Hilbert scalar
product on $A$ which induces a Hausdorff pre-Hilbert structure on
the quotient $D_{\ph}=A/I_{\ph}$ where $I_{\ph}=\{f\in A\vert
\ph(f^\ast f)=0\}$. In view of the Schwarz inequality, $I_{\ph}$ is
a left ideal of $A$ so one has a homomorphism of unital algebras
$\pi_{\ph}$ of $A$ into the endomorphisms of $D_{\ph}$ which is in
fact a hermitian representation of $A$ in the Hilbert space
$\Cal H_{\ph}$ obtained by completion of $D_{\ph}$ with domain
$D(\pi_{\ph})=D_{\ph}$. Let $\Om_{\ph}\in D_{\ph}$ be the
canonical image of the unit $1\in A$ under the projection
$A\rightarrow D_{\ph}=A/I_{\ph}$. Then one has
$\ph(f)=(\Om_{\ph},\pi_{\ph}(f)\Om_{\ph})$ for any $f\in A$
and $D_{\ph}=\pi_{\ph}(A)\Om_{\ph}$. This construction which
associates to a positive linear form $\ph$ on $A$ the triplet
$(\pi_{\ph},\Cal H_{\ph},\Om_{\ph})$ of a hermitian
representation $\pi_{\ph}$ of $A$ in Hilbert space $\Cal H_{\ph}$
with $\Om_{\ph}$ in the domain of $\pi_{\ph}$ such that
$\pi_{\ph}(A)\Om_{\ph}$ is dense in $\Cal H_{\ph}$ and
$\ph=(\Om_{\ph},\pi_{\ph}(\cdot)\Om_{\ph})$ is known as the
GNS construction; given $\ph$, the triplet
$(\pi_{\ph},\Cal H_{\ph},\Om_{\ph})$ is unique up to a unitary.

Given a hermitian representation $\pi$ of a $\ast$-algebra $A$ with
domain $D(\pi)$, to each vector $\Ph\in D(\pi)$ coresponds the
positive linear form $\ph$ on $A$ defined by
$\ph(f)=(\Ph,\pi(f)\Ph)$. Conversely, the GNS construction shows
that any positive linear form on $A$ can be realized in this manner.
To the action $(f,\Ph)\mapsto \pi(f)\Ph$ of $A$ on $D(\pi)$
corresponds the action $(f,\ph)\mapsto \ph_{f}$ of $A$ on the
(strict) convex cone $A^\ast_{+}$ of its positive linear forms where
$\ph_{f}$ is defined by $\ph_{f}(g)=\ph(f^\ast g f)$ for $f,g\in
A$.

\proclaim{\nmb.{1.2}. Proposition}
The following conditions \therosteritem{i} and \therosteritem{ii} 
are equivalent for a locally convex $\ast$-algebra $A$.
\roster
\item"(i)" $A$ is a $\ast$-algebra of unbounded operators in Hilbert 
       space $\Cal H$ with domain $D$ and its locally convex topology 
       is generated by seminorms $f\mapsto \parallel f\Ph \parallel$, 
       $\Ph\in D$. 
\item"(ii)" There is a subset $\Cal S$ of positive linear forms on 
       $A$ which is invariant by the action of $A$ on $A^\ast_{+}$ 
       and which is such that the locally convex topology of $A$ is 
       generated by the seminorms $f\mapsto (\ph(f^\ast f))^{1/2}$, 
       $\ph\in \Cal S$ and is Hausdorff. 
\endroster
\endproclaim

\demo{Proof}
$(i) \Rightarrow (ii)$. This is obvious by
taking $\Cal S=\{f\mapsto (\Ph,f\Ph)\vert \Ph\in D\}$. 

$(ii)\Rightarrow (i)$. Let
$(\pi_{\ph},\Cal H_{\ph},\Om_{\ph})$ denote the GNS triplet
associated to $\ph\in \Cal S$. Take $\Cal H$ to be the Hilbertian
direct sum $\hat\oplus_{\ph\in \Cal S}\Cal H_{\ph}$, take
$D=\oplus_{\ph\in \Cal S}\pi_{\ph}(A)\Om_{\ph}$ and notice
that it follows from the assumptions that
$\pi=\oplus_{\ph\in\Cal S}\pi_{\ph}$ is injective so $A$
identifies canonically to the $\ast$-algebra $\pi(A)$ of unbounded
operators in $\Cal H$ with domain $D$. It is clear that the locally
convex topology on $A$ generated by the seminorms $f\mapsto
(\ph(f^\ast f))^{1/2}$, $\ph\in \Cal S$ is the same as the one
generated by the seminorms $f\mapsto \parallel \pi(f)\Ph\parallel$,
$\Ph\in D$.
\qed\enddemo

Notice that if $\ph$ is a positive linear form on $A$ one has
$$
\vert\ph(f)\vert\leq (\ph(1))^{1/2}(\ph(f^\ast f))^{1/2}
$$
for any $f\in A$ (Schwarz inequality) so any $\ph\in\Cal S$ is
automatically continuous, (notice also that the same inequality shows
that $\ph=0$ whenever $\ph(1)=0$ for $\ph\in A^\ast_{+}$).

\subhead\nmb.{1.3}. Definition \endsubhead
Let $A$ be a $\ast$-algebra, $\Cal S$ be a subset of positive linear
forms on $A$ invariant by the action of $A$ on $A^\ast_{+}$ and let
$\Cal D$ be a Lie subalgebra of the Lie algebra $\Der(A)$ of
derivations of $A$ which is also a $Z(A)$-submodule of $\Der(A)$
where $Z(A)$ denotes the center of $A$. Assume that :
\roster
\item The locally convex topology on $A$ generated by the semi-norms 
       $f\mapsto \nu_{\ph}(f)=(\ph(f^\ast f))^{1/2}$, $\ph\in\Cal S$ 
       is Hausdorff;
\item $\cap\{\ker(X)\vert X\in \Cal D\}=\Bbb C1$;
\item The locally convex topology $\ta(\Cal S,\Cal D)$ on $A$ 
       generated by the seminorms 
       $\nu_{\ph}\circ X_{1}\circ \dots X_{p}$, $\ph\in\Cal S$, 
       $X_{i}\in\Cal D$, $p\in\Bbb N$ is such that 
       $(A,\ta(\Cal S,\Cal D))$ is complete. 
\endroster
Then $A$ will be said to be a {\it smooth $\ast$-algebra} relative
to $\Cal S$ and $\Cal D$, or simply a smooth $\ast$-algebra when no
confusion arises, the topology $\ta=\ta(\Cal S,\Cal D)$ being called
{\it smooth topology} of $A$.

\subhead\nmb.{1.4}. Commutative smooth $*$-algebras \endsubhead
Let $M$ be a smooth finite dimensional manifold, let 
$A=C^\infty(M,\Bbb C)$ be the $*$-algebra of all complex valued 
smooth functions on $M$. Let $\Cal D=\Der(A)=\X(M)\otimes \Bbb C$ be 
the Lie algebra of all derivations of $C^\infty(M,\Bbb C)$, i.e\. all 
complex valued vector fields on $M$. Let $\on{Vol}(M)\to M$ be 
the real line bundle of all densities on $M$, and let 
$\Ga_c^+(\on{Vol}(M))$ be the space of all smooth non-negative 
densities with compact support on $M$. Let $\Cal S$ be the space 
of all linear functionals of the form 
$f\mapsto \int_M f\mu$ for all $\mu\in \Ga_c^+(\on{Vol}(M))$. 

Then the locally convex topology on $C^\infty(M,\Bbb C)$ described in 
\nmb!{1.3.1} is the compact open topology which is Hausdorff. 
Condition \nmb!{1.3.2} is obviously satisfied. The topology 
$\ta(\Cal S,\Cal D)$ from \nmb!{1.3.3} is equivalent to the compact 
$C^\infty$-topology, i.e\. the topology of uniform convergence on 
compact subsets in all derivatives separately. Thus 
$C^\infty(M,\Bbb C)$ is a smooth $*$-algebra relative to $\Cal S$ and 
$\Cal D$.  

\head\totoc\nmb0{2}. The non-commutative torus \endhead

\subhead\nmb.{2.1}. The non-commutative torus \endsubhead
By Fourier expansion the algebra $C^\infty(S^1\x S^1,\Bbb C)$ of all 
smooth functions on the torus consists of all 
$$
f=\sum_{(k,l)\in\Bbb Z\x\Bbb Z} f_{k,l}u^kv^l,
\tag{\nmb:{1}}$$
where $(f_{k,l})$ is any rapidly decreasing sequence of complex 
numbers, i\.e\. for each $m\in \Bbb N$ the seminorm
$$
\|f\|_m:=\sup_{k,l\in\Bbb Z}{|f_{k,l}|}\,{(1+|k|+|l|)^{m}}
     < \infty,\tag{\nmb:{2}}
$$
and where $u=\exp(2\pi i t)$ and $v=\exp(2\pi i s)$ are the 
coordinates on the torus. 

Let us fix a complex number $q$ with $|q|=1$. Then the smooth 
$q$-torus $C^\infty(T^2_q)$ is the convenient associative algebra (in 
fact a Fr\'echet algebra) which is given by all elements of  the form 
\thetag{\nmb|{1}}, but where we assume now that $U$, $V$  are two 
indeterminates which satisfy
$$
UV=qVU\tag{\nmb:{3}}.
$$
Defining 
$$
U^*:=U\i,\quad V^*:=V\i
\tag{\nmb:{4}}
$$
makes $C^\infty(T^2_q)$ into a $*$-algebra. 
Note that 
$U^kV^l=q^{kl}V^lU^k$
and hence 
$$\align
fg &= \Bigl(\sum_{k,l} f_{k,l}U^kV^l\Bigr)\Bigl(\sum_{m,n} 
     g_{m,n}U^mV^n\Bigr)
     = \sum_{k,l}\Bigl(\sum_{m,n} f_{m,n}g_{k-m,l-n}q^{-n(k-m)} 
     \Bigr)U^kV^l \\
f^* &= \Bigl(\sum_{k,l} f_{k,l}U^kV^l\Bigr)^* 
     = \sum_{k,l} \bar f_{k,l}V^{-l}U^{-k}
     = \sum_{k,l} \bar f_{-k,-l}q^{-kl}U^kV^l.
\endalign$$
Using the convention
$$
f =\sum_{k,l} f'_{k,l}q^{-\frac{kl}2}U^kV^l,\quad 
     \text{ so }f'_{k,l}=f_{k,l}q^{\frac{kl}2}
$$
we get nicer descriptions for the product and the adjoint $f^*$:
$$\align
fg &= \Bigl(\sum_{k,l} f'_{k,l}q^{-\frac{kl}2}U^kV^l\Bigr)
     \Bigl(\sum_{k,l}g'_{m,n}q^{-\frac{mn}2}U^mV^n\Bigr)  \\
&= \sum_{k,l}\Bigl(\sum_{m,n} f'_{m,n}g'_{k-m,l-n}q^{-\frac{1}2(kn-ml)}
     \Bigr)q^{-\frac{kl}2}U^kV^l \\
f^* &= \sum_{k,l} \bar f'_{-k,-l}q^{-kl/2}U^kV^l.
\endalign$$

If (the argument of) $q$ is rational (mod $2\pi$), let 
$N\in\Bbb N$ be the smallest positive natural number such that 
$q^N=1$. If $q$ is irrational, we put $N=0$.

\proclaim{\nmb.{2.2}. Proposition}
If $q$ is rational, then there exits a smooth vector bundle 
$A_q\to S^1\x S^1$ with standard fiber the algebra 
$\operatorname{Mat}_N(\Bbb C)$ of all complex $(N\x N)$-matrices and 
with transition functions in $GL(n,\Bbb C)$ acting on 
$\operatorname{Mat}_N$ by conjugation, such that the non-commutative 
torus $C^\infty(T^2_q)$ is isomorphic to the algebra $\Ga(A_q)$ of all 
smooth sections of the algebra bundle $A_q\to S^1\x S^1$. The center 
of $C^\infty(T^2_q)$ is isomorphic to $C^\infty(S^1\x S^1,\Bbb C)$.
The first Chern class of the complex vector bundle $A_q$ vanishes.

Moreover, there is a smooth vector bundle $E_q\to S^1\x S^1$ with 
standard fiber $\Bbb C^N$ such that $A_q$ is the full endomorphism 
bundle $\operatorname{End}(E_q)$. The first Chern class of $E_q$ 
also vanishes.
\endproclaim

\demo{Proof}
We first claim that the algebra $\operatorname{Mat}_N$ is the unique 
algebra generated by two unitary elements $U_0$ and $V_0$ which are 
subject to the relations
$$
U_0.V_0=qV_0.U_0,\quad U_0^N=V_0^N = \Bbb I.
\tag{\nmb:{1}}$$
To see this note that each element in the algebra generated by $U_0$ 
and $V_0$ may be written in the form 
$\sum_{0\le k,l\le N-1} a_{k,l}U_0^kV_0^l$, so this algebra is of 
dimension $\le N^2$. On the other hand we consider the matrices in 
$\operatorname{Mat}_N$,
$$
U_0 = \pmatrix 0 & 1 & 0 & 0 & \dots  \\ 
               0 & 0 & 1 & 0 &\dots  \\
               \vdots& &\ddots&\ddots \\
               0 &\dots&  & 0 &  1\\
               1 & 0 &\dots&  & 0\endpmatrix,\quad
V_0 = \pmatrix 1 & 0 & 0 & \dots  \\ 
               0 & q & 0 & \dots  \\
               0 & 0 & q^2 & 0 & \dots  \\
               \vdots  & & &\ddots& \\
               0&\dots & & 0 & q^{N-1}\endpmatrix, 
$$
which satisfy relations \thetag{\nmb|{1}} and thus generate a 
$C^*$-subalgebra which  
clearly commutes only with the multiples of the identity, so it has 
to be the full matrix algebra.

Now we consider the trivial bundle
$S^1\x S^1 \x \operatorname{Mat}_N @>{pr_{1,2}}>> S^1\x S^1$.
The space of smooth section is then 
$C^\infty(S^1\x S^1,\operatorname{Mat}_N) = 
C^\infty(S^1\x S^1,\Bbb C)\otimes \operatorname{Mat}_N$,
which is generated by the unitary central elements $u$, $v$, and 
unitary $U_0$, $V_0$ with the relations \thetag{\nmb|{1}}, where the 
coefficients are again rapidly decreasing with respect to the powers 
of $u$ and $v$.

Consider now the cyclic group $\Bbb Z_N=\Bbb Z/N.\Bbb Z$, the 
$q$-action of $(m,n)\in \Bbb Z_N\x \Bbb Z_N = \Bbb Z_N^2$ on $S^1\x S^1$ 
given by $(u,v)\mapsto (q^m.u,q^n.v)$, and the $q$-action on 
$\operatorname{Mat}_N$ given by 
$A\mapsto U_0^n.V_0^{-m}.A.V_0^{m}.U_0^{-n}$. 
Note that inside the adjoint action of $GL(n,\Bbb C)$ the matrices 
$U_0$ and $V_0$ commute, since they do so in $PGL(n,\Bbb C)$, and 
that $(m,n)$ maps $U_0$ to $q^m.U_0$, and maps $V_0$ to $q^n.V_0$.
We may consider the 
following diagram, where the horizontal arrows are covering mappings
since all involved actions are strictly discontinuous, 
and where the left vertical arrow is $\Bbb Z_N^2$-equivariant.
$$\CD
S^1\x S^1 \x \operatorname{Mat}_N @>{\Bbb Z_N^2}>> A_q \\
@V{pr_{1,2}}VV  @VV{p_q}V\\
S^1\x S^1 @>{\pi}>{\Bbb Z_N^2}> S^1\x S^1
\endCD$$
Since the action of $\Bbb Z_N^2$ on $\operatorname{Mat}_N$ is by 
algebra automorphisms, the resulting smooth mapping 
$A_q\to S^1\x S^1$ is a smooth algebra bundle. The sections of $A_q$ 
correspond exactly to the $\Bbb Z_N^2$-equivariant sections of the 
left hand side. 
A section $f:S^1\x S^1\to \operatorname{Mat}_N$, 
$$
f=\sum \Sb k,l\in\Bbb Z\\ 0\le s,t\le N-1 \endSb 
     c_{k,l,s,t}\,u^kv^lU_0^sV_0^t
$$
is $\Bbb Z_N^2$-equivariant if and only if the following condition is 
satisfied:
\roster
\item "" $c_{k,l,s,t}\ne 0$ only if $k \equiv s \mod N$ and 
       $l \equiv t \mod N$.
\endroster
But then we may put $c_{k,l}=c_{k,l,s,t}$, where $s\equiv k\mod N$ and $t\equiv l\mod N$, 
and the 
section $f$ can be written as
$$f=\sum_{k,l\in\Bbb Z} c_{k,l}(uU_0)^k(vV_0)^l.$$
We just have to note that $U=uU_0$ and $V=vV_0$ satisfy only the relations 
\nmb!{2.1.3} of the noncommutative torus.

The first Chern class $c_1(A_q)$ of the complex vector bundle $A_q$ 
vanishes, by the following argument:
The mapping $\pi:S^1\x S^1\to S^1\x S^1$ in the diagram above is an 
$N^2$-sheeted covering, has mapping degree $N^2$. Thus the 
mapping in cohomology is 
$H^2(\pi)=N^2:H^2(S^1\x S^1,\Bbb Z)=\Bbb Z\to H^2(S^1\x S^1,\Bbb Z)=\Bbb Z$.
We have $H^2(\pi)c_1(A_q) = c_1(\pi^*A_q) = 0$ since $\pi^*A_q$ is a 
trivial bundle. Thus also $c_1(A_q)=0$.

Now we will construct the bundle $E_q\to S^1\x S^1$. We cannot push 
it down from a trivial bundle via the group action by $\Bbb Z_N^2$ 
since $w\mapsto U_0^n.V_0^{-m}.w$ is not a representation of 
$\Bbb Z_N^2$ on $\Bbb C^N$. We have to absorb the non-commutativity 
into a larger group acion. Thus we consider the following semidirect 
product group, its action on $S^1\x S^1\x S^1$, and its unitary 
representation on $\Bbb C$:
$$\gather
S^1\to (\Bbb Z_N \x \Bbb Z_N) \ltimes S^1\to \Bbb Z_N \x \Bbb Z_N\\
(m,n,\th).(m',n',\th') = (m+m', n+n',\th\th'q^{mn'}) \\
((\Bbb Z_N \x \Bbb Z_N) \ltimes S^1)\x(S^1\x S^1\x S^1) 
     \to S^1\x S^1\x S^1 \\
(m,n,\th).(\ph,\ps,\nu) = (q^m\ph,q^n\ps,\th\nu\ps^m)\\
(m,n,\th).w = \th U_0^n V_0^{-m}.w,\quad w\in \Bbb C^N.
\endgather$$
Using the actions we can define the bundle $E_q\to S^1\x S^1$ as follows:
$$\CD
S^1\x S^1\x S^1 \x \Bbb C^N @>{\Bbb Z_N \x \Bbb Z_N \x S^1}>> E_q \\
@V{pr_{1,2}}VV  @VV{p_q}V\\
S^1\x S^1\x S^1 @>{\Bbb Z_N \x \Bbb Z_N \x S^1}>> S^1\x S^1
\endCD$$
It is easy to check that all these actions are compatible with each other 
in such a way that we get a free fiberwise action of the algebra bundle 
$A_q$ on the vector bundle $E_q$. By counting dimensions we see that 
$A_q=\operatorname{End}(E_q)$.
For the first Chern class we can repeat the argument from above. 
\qed\enddemo

\proclaim{\nmb.{2.3}. Corollary} 
Let $q$ be a primitive $N$-th root of unity. Then the noncommutative 
torus algebra $C^\infty(T^2_q)$ is Morita equivalent to the 
commutative torus algebra $C^\infty(T^2)$.
\endproclaim

\demo{Proof} 
By theorem \nmb!{2.2} we have the algebra isomorphism 
$C^\infty(T^2_q)\cong \Ga(\operatorname{End}(E_q))$. But for any 
vector bundle the full automorphism algebra, which acts from the left 
on the space of smooth sections of the vector bundle, 
is Morita equivalent to the algebra of smooth functions on the base, 
which we may view as acting from the right. 
\qed\enddemo

\subhead\nmb.{2.4}. Derivations of the non-commutative torus \endsubhead
Let $D\in\Der(C^\infty(T^2_q))$, let us assume that $D$ is bounded. 

Then $D$ is uniquely determined by the values 
$$
D(U)=\sum_{k,l}u_{k,l}U^kV^l,\quad D(V)= 
     \sum_{k,l}v_{k,l}U^kV^l.
\tag{\nmb:{1}}$$
The relation $D(U).V+U.D(V)=qD(V).U+qV.D(U)$, by comparison of 
coefficients, leads quickly to
$$
u_{k,l-1}(1-q^{1-k}) + v_{k-1,l}(1-q^{1-l})=0.
\tag{\nmb:{2}}$$
Now let $N$ be the smallest integer with $q^N=1$ for rational $q$, 
let $N=0$ for irrational $q$.

Then for $k\equiv 1 (\mod N)$ equation \thetag{\nmb|{2}}
implies that we have $v_{k-1,l}=0$ for
$l\not\equiv 1 \mod N$, and that $v_{k-1,l}$ can be prescribed 
arbitrarily (but rapidly decreasing) for $l \equiv 1 \mod N$. 
This means that we may prescribe $D(V)=g(U^N,V^N)V$ for arbitrary
$g\in C^\infty(S^1\x S^1,\Bbb R)$.

Similarly for $l\equiv 1 \mod N$ equation \thetag{\nmb|{2}} 
implies that we have $u_{k,l-1}=0$ for 
$k\not\equiv1\mod N$, and that $u_{k,l-1}$ can be prescribed 
arbitrarily (but rapidly decreasing) for $k \equiv 1 (\mod N)$. 
This means that we may prescribe $D(U)=f(U^N,V^N)U$ for arbitrary 
$f\in C^\infty(S^1\x S^1,\Bbb R)$.

Let us write $D_U$ for the derivation given by $D_U(U)=U$ and 
$D_U(V)=0$, similarly $D_V\in\Der(C^\infty(T^2_q))$ is given by $D_V(U)=0$ 
and $D_V(V)=V$.
Thus for any $f, g \in C^\infty(S^1\x S^1,\Bbb C)$, in the 
center of $C^\infty(T^2_q)$, the expression 
$$
f(U^N,V^N)D_U + g(U^N,V^N)D_V
\tag{\nmb:{3}}$$
describes a derivation which is not inner, since it acts on the 
center (if $N>0$).

On the other hand for any $a=\sum a_{k,l}U^kV^l$ the inner 
derivation $\ad(a)b=a.b-b.a$ satisfies
$$\align
\ad(a)U &= \sum_{k,l} a_{k-1,l}(q^{-l}-1)U^kV^l,\\
\ad(a)V &= \sum_{k,l} a_{k,l-1}(1-q^{-k})U^kV^l,\\
\endalign$$
so that all other derivations specified by \thetag{\nmb|{2}} are inner 
derivations.

So we see that $\Der(C^\infty(T^2_q))= \Int(C^\infty(T^2_q)) 
\rtimes \Out(C^\infty(T^2_q))$, a semidirect product with 
$\Int(C^\infty(T^2_q))$ an ideal, where the action of 
$\Out(C^\infty(T^2_q))$ on $\Int(C^\infty(T^2_q))$ (the same as on 
$C^\infty(T^2_q)$) is given by the expression \thetag{\nmb|{2}}.

For $q$ rational, the description \thetag{\nmb|{3}} corresponds to the 
covariant derivative $\nabla_X$ along the vectorfield 
$X=f(u,v)\tfrac{\partial}{\partial t} + 
     g(u,v)\tfrac{\partial}{\partial s}$
on $S^1\x S^1$, where $u=e^{2\pi i t}$ and $v=e^{2\pi i s}$, with 
respect to the unique flat connection on the algebra bundle 
$A_q\to S^1\x S^1$, which is induced by the description in 
\nmb!{2.2}, and which respects the fiberwise `matrix'-multiplication. 
In this case the outer derivations correspond exactly to the 
derivations of the center. 

For $q$ irrational this is not the case. Here $\Out(C^\infty(T^2_q))$ is 
linearly generated by the two derivations $D_U$ and $D_V$.

\subhead\nmb.{2.5}. Conjecture \endsubhead
It might be the case that every (algebraic) derivation of the 
non-commutative torus is automatically bounded. This would follow 
from an automatic continuity result for algebra homomorphisms. 
One can find such results in the literature but they have too strong 
assumptions to be immediately applicable.  

The following argument shows how to carry over continuity from algebra 
homomorphisms to derivations:
A linear mapping $D:C^\infty(T^2_q)\to C^\infty(T^2_q)$ is a 
derivation if and only if the mapping 
$(\Id,D\ep):C^\infty(T^2_q)\to C^\infty(T^2_q)\x C^\infty(T^2_q)\ep$ is 
an algebra homomorphism, where $\ep$ is in the center and 
$\ep^2=0$ so that the multiplication 
in $C^\infty(T^2_q)\x C^\infty(T^2_q)\ep$ is given by 
$(f+g\ep)(f'+g'\ep)=ff'+(fg'+gf')\ep$.

\subhead\nmb.{2.6}. The non-commutative torus is a smooth 
$*$-algebra \endsubhead 
In fact we will show that the topology described in \nmb!{1.1.3} is 
the one we started with in \nmb!{2.1}. 

What are the states on $C^\infty(T^2_q)$? 
We consider first the trace 
$\operatorname{tr}(\sum_{k,l} c_{k,l}U^kV^l)= c_{0,0}$.
We will use only states of the form 
$$
f\mapsto \om_g(f)=\operatorname{tr}(g^*fg)
$$
for some $g\in C^\infty(T^2_q)$, and indeed $g=1$ will suffice.
We start  
to check that we can reproduce a generating system of seminorms. For 
that it suffices to consider 
$$\multline
f=\sum_{k,l} c_{k,l}U^kV^l \mapsto \om_1(f)=
\operatorname{tr}(f^*f)^{1/2} = \\
=\operatorname{tr}\Bigl(\sum_{k,l,m,n}  
 \overline{c_{m,n}}V^{-n}U^{-m}c_{k,l}U^kV^l\Bigr)^{1/2}
= \Bigl(\sum_{k,l} \overline{c_{k,l}} c_{k,l}\Bigr)^{1/2} = 
\|f\|_{\ell^2}
\endmultline$$
and to compose it with an appropriate composition of the two basic 
derivations $D_U$ and $D_V$ from \nmb!{2.4} which give us:
$$
D_U^mD_V^n\Bigl(\sum_{k,l} c_{k,l}U^kV^l\Bigr) = \sum_{k,l} 
c_{k,l}k^ml^nU^kV^l. 
$$

It remains to show that an arbitrary state $\om$ on $C^\infty(T^2_q)$ is 
bounded: We use the Gelfand-Naimark-Segal construction. The subspace 
$I_\om:=\{f\in C^\infty(T^2_q):\om(f^*f)=0\}$ is a left ideal, since 
by the Cauchy-Schwarz inequality we have 
$\om((gf)^*gf)=\om((f^*g^*g)f)\le \om(g^*gff^*g^*g)\om(f^*f)=0$.
Then $D_\om:=C^\infty(T^2_q)/I_\om$ is a pre-Hilbert space with 
the inner product $\om(f^*g)$ which is positively defined by the 
definition of $I_\om$. We get a *-representation  
$\pi_\om:C^\infty(T^2_q)\to L(D_\om,D_\om)$. Since 
$\Id_{D_\om}=\pi_\om(U^*U)=\pi_\om(U)^*\pi_\om(U)=\pi_\om(U)\pi_\om(U)^*$,
the operators $\pi_\om(U)$ and $\pi_\om(V)$ are unitary. Since 
the coefficients in $C^\infty(T^2_q)$ are rapidly decreasing, 
$$
\pi_\om(f)=\pi_\om(\sum_{k,l} c_{k,l}U^kV^l)
     =\sum_{k,l} c_{k,l}\pi_\om(U)^k\pi_\om(V)^l
$$
is a bounded operator for each $f\in C^\infty(T^2_q)$, and $\pi_\om$ 
is bounded. 
Thus the representation $\pi_\om$ and the state $\om$ can be 
extended to the `$C^*$-algebra completion'  
$C(T^2_q)$ of $C^\infty(T^2_q)$ and $\om$ has norm 1 on 
$C(T^2_q)$. Since $C^\infty(T^2_q)\to C(T^2_q)$ is continuous, $\om$ 
is bounded on $C^\infty(T^2_q)$.

\subhead
\nmb.{2.7}. Higher dimensional non-commutative tori 
\endsubhead
Let us fix a complex number $q$ with $|q|=1$, and
let us consider the algebra $C^\infty(T^n_q)$ consisting of all 
$$
f=\sum_{k=(k_1,\dots,k_n)\in \Bbb Z^n}f_{k}
     S_1^{k_1}S_2^{k_2}\dots S_n^{k_n}
\tag{\nmb:{1}}$$
where $(f_k)$ is any rapidly decreasing sequence of complex numbers
so that for each $m\in \Bbb N$ the seminorm 
$$
\|f\|_m:=\sup_{k=(k_1,\dots,k_n)\in \Bbb Z^n}
     |f_{k}|(1+|k_1|+\dots+|k_n|)^m <\infty,
$$
and where the generators $S_1,\dots,S_n$ satisfy the commutation 
rules
$$
\cases S_iS_{i+1} = qS_{i+1}S_i & \text{ for }i=1,\dots,n-1\\
       S_iS_j=S_jS_i            & \text{ for }|i-j|\ge 2 
\endcases
$$
This looks like an interesting generalization of the non-commutatve 
torus $C^\infty(T^2_q)$. But it is not so interesting as the following result 
shows:

\proclaim{Lemma}
For $n=2p$ we have 
$$
C^\infty(T^{2p}_q) = C^\infty(T^2_q)\hat\otimes \dots\hat\otimes 
     C^\infty(T^2_q) \quad p\text{ times,}
$$
where we may use the projective tensor product.

For $n=2p+1$ we have 
$$
C^\infty(T^{2p+1}_q) 
     = C^\infty(T^2_q)\hat\otimes \dots\hat\otimes C^\infty(T^2_q) 
     \hat \otimes C^\infty(S^1), 
$$
the projective tensor product of 2p copies of the non-commutative 
2-torus with one 1-torus.
\endproclaim

\demo{Proof}
Let first $n=2p$. Consider the new set of generators of the algebra 
$T^{2p}_q$
$$
\cases U_j := S_1S_3\dots S_{2j-1} &\quad \text{ for }j=1,\dots ,p\\
       V_j := S_{2j}               &\quad \text{ for }j=1,\dots ,p
\endcases
$$
Then obviously $U_jV_j=qV_jU_j$ and all other pairs commute so that 
the first result follows.

If we have moreover an element $S_{2p+1}$ then we also consider the last 
generator $Z=S_1S_3\dots S_{2p+1}$ which lies in the center of 
$T^{2p+1}_q$ (it even generates the center if $q$ is irrational) 
and thus splits off a central subalgebra isomorphic to 
$C^\infty(S^1)$.
\qed\enddemo

\subhead\nmb.{2.8}. The non-commutative circle \endsubhead
We look for the non-commutative circle as a smooth algebra which 
is a quotient of the non-commutative torus. 
Since $C^\infty(T^2_q)$ is a simple 
algebra for irrational $q$ we will succeed only for rational $q$, 
thus let us take $q\in S^1\subset \Bbb C$ with $q^N=1$ for minimal 
$N$. As in \nmb!{2.1} let $u=\exp(2\pi i t)$ and $v=\exp(2\pi i s)$ 
be the coordinates on the torus $S^1\x S^1$, and let $z=\exp(2\pi i 
x)$ be the coordinate on $S^1$. Let us consider the embedding 
$$
i:S^1\to S^1\x S^1, \quad i(z)=(z^a,z^b),
$$
where $a,b\in \Bbb Z$ are relatively prime. Then we consider the 
algebra bundle $A_q\to S^1\x S^1$ with typical fiber 
$\operatorname{Mat}_N$ constructed in the proof of proposition 
\nmb!{2.2}, and take the pullback bundle $i^*A_q\to S^1$ and the 
space of smooth sections is then viewed as the 
\idx{\it non-commutative $q$-circle}. We want to describe it by 
generators and relations. For that consider the following diagram 
$$\cgaps{;}\newCD
S^1\x S^1\x \operatorname{Mat}_N 
     @()\l{\Bbb Z_N^2}@(1,-1) @()\L{pr_{1,2}}@(0,-3) & & & 
     S^1\x \operatorname{Mat}_N @()\l{\Bbb Z_N}@(-1,-1) @()\L{}@(-3,0) 
     @()\l{pr_1}@(0,-3) \\
& A_q @()\L{p}@(0,-1) & i^*A_q @()\L{p^*i}@(-1,0) @()\L{i^*p}@(0,-1) & \\
& S^1\x S^1 & S^1 @()\L{i}@(-1,0) & \\
S^1\x S^1 @()\L{\Bbb Z_N^2}@(1,1) & & & 
     S^1 @()\L{i}@(-3,0) @()\L{\Bbb Z_N}@(-1,1) \\
\endnewCD$$
where all diagonal mappings are covering maps with the groups of 
covering transformations indicated: $(m,n)\in\Bbb Z_N^2$ 
acts on $S^1\x S^1$ by $(u,v)\mapsto(q^mu,q^nv)$ and on 
$\operatorname{Mat}_N$ by $A\mapsto U_0^n.V_0^{-m}.A.V_0^m.U_0^{-n}$;
$p\in\Bbb Z_N$ acts on $S^1$ by $z\mapsto q^pz$ and on 
$\operatorname{Mat}_N$ by 
$A\mapsto U_0^{bp}.V_0^{-ap}.A.V_0^{ap}.U_0^{-bp}$.
The outer horizontal mappings are equivariant with respect to the 
homomorphism $\Bbb Z_N\to \Bbb Z_N^2$ which is given by 
$p\mapsto (ap,bp)$. So the smooth sections of the algebra bundle 
$i^*A_q\to S^1$ correspond to the $\Bbb Z_N$-equivariant smooth 
functions $S^1\to \operatorname{Mat}_N$. A smooth function  
$$
f=\sum \Sb k\in\Bbb Z\\ 0\le s,t\le N-1 \endSb 
     c_{k,s,t}\,z^kU_0^sV_0^t
$$
is $\Bbb Z_N$-equivariant if and only if the following condition is 
satisfied:
$$
c_{k,s,t}\ne 0 \text{  only if }k \equiv as+bt\mod N
$$ 
But then the 
function $f$ can be written as
$$\align
f &=\sum \Sb j\in\Bbb Z\\ 0\le s,t\le N-1 \endSb  
     c_{as+bt+jN,s,t}z^{jN}(z^aU_0)^s(z^bV_0)^t\\
&=\sum \Sb j\in\Bbb Z\\ 0\le s,t\le N-1 \endSb  
     c_{as+bt+jN,s,t}Z^jU^sV^t,\\
\endalign$$
where 
$Z:=z^N$, $U:=z^aU_0$ and $V:=z^bV_0$ satisfy the relations 
$$
UV=qVU,\quad Z \text{ is central,}\quad U^{N}=Z^a,\quad V^{N}=Z^b.
\tag{\nmb:{1}}$$
We also have $Z=U^{Na'}V^{Nb'}$ where $a',b'\in\Bbb Z$ satisfy $aa'+bb'=1$.
So the \idx{\it non-commutative $q$-circle of slope $b/a$} 
in the non-commutative $q$-torus is the associative algebra 
generated by two elemets $U,V$ with the relations \thetag{\nmb|{1}}, 
and with  
rapidly decreasing coefficients. 

If $q=1$ we have $N=1$, thus $U=Z^a$, $V=Z^b$, and 
clearly we just have the algebra of smooth functions on $S^1$.

\head\totoc\nmb0{3}. The smooth Heisenberg algebra \endhead

\subhead\nmb.{3.1} \endsubhead
We recall here (see \cit!{22}, \cit!{29}, or \cit!{30}) 
some wellknown results from the theory of 
distributions which we shall need in the following.
We consider the following spaces of smooth functions on $\Bbb R^n$:

The space $\Cal S(\Bbb R^n)$ of all rapidly decreasing smooth 
functions $f$ for which 
$x\mapsto(1+|x|^2)^k\partial^\al f(x)$ is bounded for all  
$k\in \Bbb N$ and all multiindices $\al\in \Bbb N_0^n$, with the 
locally convex topology described by these conditions, a nuclear 
Fr\'echet space.  Its dual space $\Cal S'(\Bbb R^n)$ is the space of 
tempered distributions.

The space $\Cal O_C(\Bbb R^n)$ of all smooth functions $f$ on 
$\Bbb R^n$ for which there exists $k\in \Bbb Z$ such that 
$x\mapsto (1+|x|^2)^k\partial^\al f(x)$ is bounded for each multiindex 
$\al\in\Bbb N_0^n$, with the locally convex topology described by 
this condition (a nuclear LF space). Its dual space 
$\Cal O_C'(\Bbb R^n)$ is usually called the space of rapidly 
decreasing distributions (see \cit!{29}).

The space $\Cal O_M(\Bbb R^n)$ of all smooth functions $f$ 
on $\Bbb R^n$ such that for each multiindex 
$\al\in\Bbb N_0^n$  there exists $k\in \Bbb Z$ such that 
$x\mapsto(1+|x|^2)^k\partial^\al f(x)$ is bounded, with the locally convex 
topology described by this condition (a nuclear space). This is the 
space of tempered smooth functions.
Its dual space $\Cal O_M'(\Bbb R^n)$ will be called the space of 
\idx{\it speedily decreasing distributions}.

There are the following inclusions between these spaces: 
$$
\Cal S\subset \Cal O_C\subset \Cal O_M\subset \Cal S',\qquad
\Cal S\subset \Cal O_M'\subset \Cal O_C'\subset \Cal S'.
$$ 
The Fourier transform of functions $f\in \Cal S$ and its inverse,
$$
\Cal Ff(y):= \int_{\Bbb R^n}e^{-i\langle x,y\rangle}f(x)\,dx, \quad
\Cal F\i f(x):= 
     \tfrac 1{(2\pi)^n}\int_{\Bbb R^n}e^{i\langle x,y\rangle}f(y)\,dy
$$ 
extend to isomorphisms of $\Cal S'$, which induce isomorphisms 
$\Cal O_M\to \Cal O_C'$ and $\Cal O_C\to \Cal O_M'$.
Under the convolution product $(f*g)(x)=\int_{\Bbb R^n}f(x-y)g(y)dy$ 
the space $\Cal S$ is a commutative algebra and the Fourier transform 
is an isomorphism between this and the pointwise multiplication. The 
convolution carries over to distributions as follows: It induces an 
associative commutative product on $\Cal O_C'$ and makes $\Cal S'$ 
into an $\Cal O_C'$-module. The Fourier transform is an algebra 
isomorphism $\Cal F: (\Cal O_C',*)\to (\Cal O_M,\cdot)$. The 
convolution $*$ is jointly continuous on $\Cal O_C'$. 
Moreover 
$\Cal S*\Cal S'\subset \Cal O_M$ and $\Cal S*\Cal O_C'\subset \Cal S$.
See \cit!{29}, pp\. 246ff and 268. 
The space $\Cal O_M \cong \Cal O_C'$ is a complete bornological 
nuclear locally convex vector space, and the dual  
$\Cal O_M' \cong \Cal O_C$ is a complete nuclear (LF)-space, thus 
also bornological, see  
\cit!{15},~II,~\S4,4,~th\'eor\`eme~16 (page 131).
Let us summarize the embeddings and isomorphisms in the following 
diagram:
$$\rgaps{0.6;0.6}\newCD
& \Cal O_M' @(1,0) &  \Cal O_C' @(1,-1) & \\
\Cal S @(1,1) @(1,-1) & & & \Cal S'\\
& \Cal O_C @(1,0) @()\L{\Cal F}\l{\cong}@(0,2) & 
     \Cal O_M @(1,1) @()\L{\Cal F}\l{\cong}@(0,2) & \\
\endnewCD$$

{\sl Moreover we have 
$\Cal O_C(\Bbb R^n) 
\hat\otimes \Cal O_C(\Bbb R^m)\cong \Cal O_C(\Bbb R^{n+m})$
and
$\Cal O_M(\Bbb R^n) 
\hat\otimes \Cal O_M(\Bbb R^m)\cong \Cal O_M(\Bbb R^{n+m})$, for the 
completed projective tensor product which agrees with the injective 
one. }

Since we have been unable to locate this result in the literature we 
sketch a proof: We start with $\Cal O_M$. By (\cit!{29}, p.\ 246) the 
space $\Cal O_M(\Bbb R^n)$ is the space of the multipliers in 
$L_b(\Cal S(\Bbb R^n),\Cal S(\Bbb R^n))$, with the induced topology, 
where $L_b$ denotes the space of continuous linear mappings with the 
topology of uniform convergence on bounded sets (i.e\. on compact 
sets, since $\Cal S$ is Montel), whose bornology is the same as that 
from \nmb!{4.5.1}. It is well known that 
$L_b(\Cal S(\Bbb R^n),\Cal S(\Bbb R^n))
     \cong \Cal S(\Bbb R^n)'\hat\otimes \Cal S(\Bbb R^n)$.
Thus we have the following diagrams of embeddings:
$$\cgaps{;}\newCD
\Cal O_M(\Bbb R^n)\hat\otimes\Cal O_M(\Bbb R^m)  @()\a-@(1,0)  
     @()\0(@(0,-1) &
     \Cal O_M(\Bbb R^{n+m})  @()\0(@(0,-1)\\
L_b(\Cal S(\Bbb R^n),\Cal S(\Bbb R^n))\hat\otimes 
     L_b(\Cal S(\Bbb R^m),\Cal S(\Bbb R^m))  @()\a=@(0,-1) &
     L_b(\Cal S(\Bbb R^{n+m}),\Cal S(\Bbb R^{n+m}))  @()\a=@(0,-1)\\ 
\Cal S(\Bbb R^n)'\hat\otimes\Cal S(\Bbb R^n)\hat\otimes
     \Cal S(\Bbb R^m)'\hat\otimes\Cal S(\Bbb R^m)  @()\a=@(1,0) &
     \Cal S(\Bbb R^{n+m})'\hat\otimes\Cal S(\Bbb R^{n+m})
\endnewCD$$
It remains to check that the spaces of smooth functions with compact 
support are dense in $\Cal O_M$, which is easy, and that the trace 
topology on subspaces of functions with fixed compact support is the usual 
Fr\'echet topology, so that 
$C^\infty_c(\Bbb R^n)\otimes C^\infty_c(\Bbb R^m)$ is dense in 
$\Cal O_M(\Bbb R^n)\hat\otimes\Cal O_M(\Bbb R^m)$. Thus the result 
for $\Cal O_M$ follows. 
For $\Cal O_C$ we get then the result by 
$\Cal O_C(\Bbb R^{n+m})\cong \Cal O_M(\Bbb R^{n+m})' 
\cong (\Cal O_M(\Bbb R^n)\hat\otimes\Cal O_M(\Bbb R^m))'
\cong \operatorname{Bilin}_{\text{cont}}
     (\Cal O_M(\Bbb R^n),\Cal O_M(\Bbb R^m);\Bbb R)\cong
\Cal O_M(\Bbb R^n)'\hat\otimes\Cal O_M(\Bbb R^m)'\cong
\Cal O_C(\Bbb R^n)\hat\otimes\Cal O_C(\Bbb R^m)$.

\subhead\nmb.{3.2}. The Heisenberg relation \endsubhead
Let $Q,P$ be two generators which satisfy the Heisenberg relation
$$
[Q,P]=QP-PQ=i\hbar.\tag{\nmb:{1}}
$$
We suppose that they are hermitian: $Q^*=Q$ and 
$P^*=P$, which implies that $\hbar$ should be real.

\proclaim{Lemma} Then the unitary generators $e^{iQ}$ and $e^{iP}$ 
satisfy the Weyl relation 
$$
e^{itQ}.e^{isP} = e^{-its\hbar}.e^{isP}.e^{itQ}
\text{ for }(t,s)\in \Bbb R^2\tag{\nmb:{2}}
$$
\endproclaim

\demo{Algebraic proof}
We claim that the Heisenberg relations imply that for all 
$m,n\in \Bbb N_0$ we have
$$
Q^nP^m = \sum_{k=0}^\infty 
     \binom nk \binom mk k!(i\hbar)^k P^{m-k}Q^{n-k},\tag{\nmb:{3}}
$$
which is in fact a finite sum. In the simplest cases 
\thetag{\nmb|{3}} boils down to
$QP^m=P^mQ + mi\hbar P^{m-1}$ and $Q^nP=PQ^n+ni\hbar Q^{n-1}$ which 
follow easily from \thetag{\nmb|{1}}.
 From these simple cases one may then prove \thetag{\nmb|{3}} by 
induction. Finally \thetag{\nmb|{2}} follows from \thetag{\nmb|{3}} 
by a simple power series calculation.
\qed\enddemo

\demo{Analytic proof}
Another proof of \thetag{\nmb|{2}} goes as follows. Let $Q$ and $P$ act on the 
space $\Cal S(\Bbb R)$ of all rapidly decreasing functions, by 
$(Qf)(u)=uf(u)$ and $(Pf)(u)=\frac{\hbar}{i}\partial_uf(u)$. Then the 
operators $Q$ and $P$ satisfy the Heisenberg relation \thetag{\nmb|{1}}, and 
they are selfadjoint with respect to the inner product 
$\int_{\Bbb R}\overline{f(u)}g(u)du$. It is more difficult to see that 
there are no other relations between these operators.
Let us consider the smooth 1-parameter subgroups of isomorphisms
$e^{isP}$ and $e^{itQ}$ with infinitesimal generators $iP$ and $iQ$:
$$\align
\left(e^{isP}f\right)(u) :&=
f(u+s\hbar),\tag{\nmb:{4}}\\
\left(e^{itQ}f\right)(u) :&= e^{itu}f(u),\\
\left(e^{isP}e^{itQ}f\right)(u) &= \left(e^{itz}f(z)\right)_{z=u+s\hbar}
     = e^{ist\hbar}e^{itu}f(u+s\hbar)\\
&= (e^{ist\hbar}.e^{itQ}.e^{isP}f)(u).\qed
\endalign$$
\enddemo

{\it Using the Baker-Campbell-Hausdorff
formula.}
Recall that (for finite dimensional matrices) we have 
$e^Qe^P=e^{C(Q,P)}$ where 
$$\align C(&Q,P) = P + \int_0^1 \sum_{n=0}^{\infty}
     \frac{(-1)^n}{n+1}(e^{t. \ad Q}.e^{ \ad P})^n.Q\,dt \\
&= Q + P + \tfrac12[Q,P] +\tfrac1{12}\Bigl([Q,[Q,P]]-[P,[P,Q]]\Bigr) 
     + \cdots
\endalign$$
Since we have $[Q,P]=i\hbar$, we see that 
$C(Q,P)=Q+P+\tfrac{i}2\hbar$.
Thus we may use formally new generating elements 
$$
e^{itQ}e^{isP}=e^{itQ+isP-\frac i2\hbar ts} = 
     e^{-\frac i2\hbar ts}e^{i(tQ+sP)}\tag{\nmb:{5}}
$$ 
and we see that the multiplication then will be 
$$\align
e^{i(x_1Q+y_1P)}e^{i(x_2Q+y_2P)} 
&= e^{-\frac i2\hbar (x_1y_2-x_2y_1)}
     e^{i((x_1+x_2)Q+(y_1+y_2)P)}\tag{\nmb:{6}}\\
&= e^{-\frac i2\hbar \om(x,y)}
     e^{i((x_1+x_2)Q+(y_1+y_2)P)},\\
\endalign$$
where $\om(x,y)=x_1y_2-x_2y_1$ is the symplectic form on $\Bbb R^2$.

\subhead\nmb.{3.3}. The twisted convolution in two versions \endsubhead
Let $Q,P$ be hermitian generators with $[Q,P]=i\hbar$ as in 
\nmb!{3.2}. For a rapidly decreasing distribution
$a(t,s)\in \Cal O_C'(\Bbb R^2)$ we consider the formal expression
$$
\int_{\Bbb R^2} a(t,s) e^{itQ}e^{isP}\,dt\,ds. \tag{\nmb:{1}}
$$
If we multiply two such expressions and compute (formally, but see 
below) in the 
space of endomorphisms of $\Cal S(\Bbb R)$
we get
$$\align
\int_{\Bbb R^2} &a(t,s) e^{itQ}e^{isP}\,dt\,ds.
\int_{\Bbb R^2} b(u,v) e^{iuQ}e^{ivP}\,du\,dv = \\
&=\int_{\Bbb R^2}\int_{\Bbb R^2} a(t,s)b(u,v) e^{itQ}e^{isP}
     e^{iuQ}e^{ivP}\,dt\,ds\,du\,dv = \\
&=\int_{\Bbb R^2}\int_{\Bbb R^2} a(t,s)b(u,v) 
     e^{isu\hbar}e^{i(t+u)Q}e^{i(s+v)P}\,dt\,ds\,du\,dv = \\
&=\int_{\Bbb R^2}\left(\int_{\Bbb R^2} a(t'-u,s'-v)b(u,v) 
e^{i(s'u-vu)\hbar}\,du\,dv\right) e^{it'Q}e^{is'P}\,dt'\,ds'
\endalign$$
so that we may consider the `twisted convolution' (formally, but see 
below)
$$
(a *_{\hbar} b)(t,s) = \int_{\Bbb R^2} a(t-u,s-v)b(u,v) 
e^{i(su-vu)\hbar}\,du\,dv.\tag 1
$$
For a speedily decreasing distribution
$a(t,s)\in \Cal O_M'(\Bbb R^2)$ we consider the formal expression
$$
\int_{\Bbb R^2} a(t,s) e^{i(tQ+sP)}\,dt\,ds =
\int_{\Bbb R^2} a(t,s) e^{\frac{i\hbar}2 ts}e^{itQ}e^{isP}\,dt\,ds. 
\tag{\nmb:{2}}
$$
If we multiply two such expressions and compute as above we get
$$\align
\int_{\Bbb R^2} &a(t,s) e^{i(tQ+sP)}\,dt\,ds.
     \int_{\Bbb R^2} b(u,v) e^{i(uQ+vP)}\,du\,dv = \\
&=\int_{\Bbb R^4} a(t,s) b(u,v) e^{i(tQ+sP)}
     e^{i(uQ+vP)}\,dt\,ds\,du\,dv = \\
&=\int_{\Bbb R^4} a(t,s) b(u,v) 
     e^{-\frac{i\hbar}2(tv-su)}e^{i((t+u)Q+(s+v)P)}
\,dt\,ds\,du\,dv = \\
&=\int_{\Bbb R^2}\left(\int_{\Bbb R^2} a(t'-u,s'-v)b(u,v) 
e^{-\frac{i\hbar}2(t'v-s'u)}\,du\,dv\right) 
e^{i(t'Q+s'P)}\,dt'\,ds',
\endalign$$
which motivates the `other twisted convolution' 
for speedily decreasing distributions 
$a,b\in \Cal O_M'(\Bbb R^{2n})$
$$
(a \hat*_\hbar b)(x)=
\int_{\Bbb R^{2n}} a(x-y)b(y) e^{-\frac{i\hbar}2\om(x,y)}\,dy\tag{\nmb:{3}}
$$
where $\om(x,y)=\sum_{i=1}^n(x_{2i-1}y_{2i}-y_{2i-1}x_{2i})$ is the 
symplectic form on $\Bbb R^{2n}$.

\proclaim{Theorem} The `twisted convolution' 
$$
(a *_{\hbar} b)(t,s) = \int_{\Bbb R^2} a(t-u,s-v)b(u,v) 
e^{i(su-vu)\hbar}\,du\,dv\tag{\nmb:{4}}
$$ 
is a well defined, jointly continuous, and associative product on the 
space $\Cal O_M'(\Bbb R^2)$ of speedily decreasing distributions. It is  
smooth in the variable $\hbar \in \Bbb R$.
The convenient algebra 
$(\Cal O_M')_\hbar:=(\Cal O_M',*_\hbar)$ is called the smooth 
Heisenberg plane with parameter $\hbar\in \Bbb R$.
The noncommutative torus $T^2_{e^{i\hbar}}$ with rotation parameter 
$q=e^{i\hbar}$ is a closed subalgebra with unit of 
$(\Cal O_M')_\hbar$, it corresponds to the subspace of all rapidly 
decreasing measures on $\Bbb R^2$ with support in the lattice 
$(2\pi\Bbb Z)^2$. The generalization of this to $\Bbb R^{2n}$ also 
holds.  

The `other twisted convolution' 
$$\align
(a \hat *_\hbar b)(x) 
&=\int_{\Bbb R^{2n}} a(x-y)b(y) e^{-\frac{i\hbar}2\om(x,y)}\,dy \tag{\nmb:{5}}
\endalign$$
is an associative bounded multiplication 
on the space $\Cal O_M'(\Bbb R^{2n})$ 
of speedily 
decreasing distributions,
and the 
algebras $(\Cal O_M'(\Bbb R^{2n}), *_\hbar)$ and 
$(\Cal O_M'(\Bbb R^{2n}),\hat *_\hbar)$ are 
isomorphic under the mapping 
$$a(x)\mapsto e^{-\frac{i\hbar}2 \sum_{i=1}^n x_{2i-1}x_{2i}}a(x).$$

Moreover, for both multiplications the algebras 
$\Cal O_M'(\Bbb R^{2n})$ decompose as 
(bornological or projective or injective) tensorproduct of $n$ 
commuting factors
$$
\Cal O_M'(\Bbb R^{2n}) = 
\Cal O_M'(\Bbb R^{2})\tilde\otimes\dots\tilde\otimes \Cal O_M'(\Bbb R^{2}).
$$

Formula \thetag{\nmb|{1}} defines a bounded linear 
mapping $\Cal O'_M(\Bbb R^{2})\to L(\Cal S(\Bbb R),\Cal S(\Bbb R))$ 
which is injective if $\hbar\ne0$, and is an algebra homomorphism from 
the twisted convolution \thetag{\nmb|{4}} to the composition. Likewise 
formula \thetag{\nmb|{2}} defines a bounded linear 
mapping $\Cal O'_M(\Bbb R^{2})\to L(\Cal S(\Bbb R),\Cal S(\Bbb R))$ 
which is injective if $\hbar\ne0$, and is an algebra homomorphism from 
the other twisted convolution \thetag{\nmb|{5}} to the composition. The 
analoga on $\Bbb R^{2n}$ also hold.
\endproclaim

\demo{Proof} 
We have to check that $a*_\hbar b$, given by \thetag{\nmb|{4}}, defines a 
distribution in $\Cal O_M'(\Bbb R^2)$. So let 
$g\in \Cal O_M(\Bbb R^2)$, then 
$$\align
\langle a*_\hbar b,g \rangle &= 
\int\int a(t-u,s-v)b(u,v)e^{i(su-vu)\hbar}g(t,s)\;du\,dv\,dt\,ds\\
:&= \int\int a(t,s)b(u,v)e^{isu\hbar}g(t+u,s+v)\;du\,dv\,dt\,ds,
\endalign$$
which makes sense since we shall see that 
$(t,s,u,v)\mapsto e^{isu\hbar}g(t+u,s+v)$ is an element in 
$\Cal O_M(\Bbb R^4)=\Cal O_M(\Bbb R^2)\hat\otimes\Cal O_M(\Bbb R^2)$, 
and moreover that 
$\hbar \mapsto ((t,s,u,v)\mapsto e^{isu\hbar}g(t+u,s+v))$ is a 
smooth curve $\Bbb R\to \Cal O_M(\Bbb R^4)$.
All this is a consequence of the following facts:
\roster
\item"(\nmb:{6})" $\Cal O_M(\Bbb R^4)$ is a bounded algebra for the pointwise 
       multiplication.
\item"(\nmb:{7})" For a polynomial $p:\Bbb R^n\to\Bbb R^m$ the mapping 
       $p^*:\Cal O_M(\Bbb R^m)\to \Cal O_M(\Bbb R^n)$ is bounded 
       linear. 
\item"(\nmb:{8})" $x\mapsto e^{ix}$ belongs to $\Cal O_M(\Bbb R)$.
\item"(\nmb:{9})" $\hbar \mapsto ((s,u)\mapsto e^{isu\hbar})$ is a smooth curve 
       in $\Cal O_M(\Bbb R^2)$ since obviously the mapping 
       $\Cal O_M(\Bbb R^3)\to C^\infty(\Bbb R,\Cal O_M(\Bbb R^2))$ is 
       bounded linear.
\endroster
This shows that $a*_\hbar b$ 
is a bounded (thus continuous, since $\Cal O_M$ is bornological by 
\cit!{15},~II,~\S4,4,~th\'eor\`eme~16 (page 131)) linear functional on 
$\Cal O_M(\Bbb R^4)$, and that $(a,b)\mapsto a*_\hbar b$ is 
bounded. 

It is easy to see that $*_\hbar$ is an associative product, since 
this is clear for $\hbar=0$ and for $\hbar\ne0$ we have an injective 
algebra homomorphism 
$(\Cal O_M'(\Bbb R^2),*_\hbar)\to L(\Cal S(\Bbb R),\Cal S(\Bbb R))$, 
see below. 

The statement about the noncommutative torus is clear.

The statement about the other twisted convolution follows via the 
isomorphism. The extension to $\Bbb R^{2n}$ is obvious and the 
decomposition into the tensorproduct follows from the considerations 
in \nmb!{3.1}.

Finally, on $\Bbb R^2$, the statement about the representation on 
$\Cal S(\Bbb R)$ can be proved as follows.
Using \nmb!{3.2.4}
for $f\in \Cal S(\Bbb R)$ we have 
$$\align
\left(\left(\int_{\Bbb R^2} a(t,s) e^{itQ}e^{isP}\,dt\,ds\right)f\right)(u) 
:&= \int_{\Bbb R^2} a(t,s) (e^{itQ}e^{isP}f)(u)\,dt\,ds = \\
&=\int_{\Bbb R^2} a(t,s) e^{itu}f(u+s\hbar)\,dt\,ds. 
\endalign$$
We observe that for $u\in \Bbb R$ and $f\in \Cal S(\Bbb R)$ the 
mapping $(t,s,u)\mapsto e^{itu}f(u+s\hbar)$ belongs to 
$\Cal O_M(\Bbb R)\hat\otimes \Cal S(\Bbb R^2)$, but not to 
$\Cal O_C(\Bbb R^2)\hat\otimes \Cal S(\Bbb R)$.
This follows from 
\therosteritem{\nmb|{6}}-\therosteritem{\nmb|{9}} and from the fact that 
for a polynomial $p:\Bbb R^n\to\Bbb R^m$ the mapping 
$p^*:\Cal S(\Bbb R^m)\to \Cal S(\Bbb R^n)$ is bounded linear. 
This implies the result, since the extension to $\Bbb R^{2n}$ is 
again obvious. 
\qed\enddemo

\subhead\nmb.{3.4}. Remarks \endsubhead
The twisted convolution $*_\hbar$ is not well 
defined on the classical space $\Cal O_C'\supset \Cal O_M'$ of 
rapidly decreasing distributions, since $e^{isu\hbar}g(t+u,s+v)$ is 
not in $\Cal O_C$, even if $g$ is in $\Cal O_C$, because 
$(s,u)\mapsto e^{isu\hbar}$ is not in $\Cal O_C$, see \cit!{29}, 
p\.~245.
Property \nmb!{3.3.7} is wrong for $\Cal O_C$, but it holds for 
linear mappings. 
Is it true that 
$\Cal O_M'$ is the optimal space of distributions on which 
the twisted convolution defines an algebra structure?  

The statement that $a*_\hbar b$ is smooth in $\hbar$ cannot be 
improved to real analytic $\Bbb R\to \Cal O_M'(\Bbb R^2)$ in the weak 
sense of \cit!{20}. The source of this is the fact that 
$\hbar\mapsto (x\mapsto e^{ix\hbar})$ is not real analytic 
$\Bbb R\to \Cal O_M(\Bbb R)$, even after composing with a linear 
functional: Let $f\in \Cal S(\Bbb R)\subset \Cal O_M'(\Bbb R)$ be 
such that the Fourier transform $\Cal Ff\in \Cal S(\Bbb R)$ is not 
real analytic. Then 
$$\hbar\mapsto \langle f,e^{i(\quad)\hbar} 
     \rangle =\int f(x)e^{ix\hbar}dx = (\Cal Ff)(-\hbar)$$
is not real analytic.
This is related to the fact that the Moyal $*$-product is only formal 
in $\hbar$, although there exist integral expressions in the sense of 
distributions which are smooth in $\hbar$, see \nmb!{3.5} and 
\nmb!{3.6} below. 

In \cit!{23} J. Maillard defined spaces of distributions 
$\Cal O_\hbar'(\Bbb R^2)$ as follows, depending on $\hbar$:  
$\Cal O_\hbar'(\Bbb R^2)$ 
consists of all distributions $a\in \Cal S'(\Bbb R^2)$ such that the 
formal expression from above
$$
\int_{\Bbb R^2} a(t,s) e^{i(tQ+sP)}\,dt\,ds
$$
defines a linear mapping $\Cal S(\Bbb R)\to \Cal S(\Bbb R)$, which 
then turns out to be bounded. 
 From \nmb!{3.3} it follows that 
$\Cal O_M'(\Bbb R^2) \subseteq \Cal O_{\hbar}'(\Bbb R^2)$.

So for the twisted convolution as in \nmb!{3.3} the (possibly) 
different spaces $\Cal O_{\hbar}'(\Bbb R^2)$ stabilize to (or at 
least contain) a fixed 
space $\Cal O_M'(\Bbb R^2)$. 
Also Kammerer in \cit!{18} gives many results on the space 
$\Cal O_{\hbar}'(\Bbb R^2)$.

\subhead\nmb.{3.5}. The Fourier transform of the twisted convolution 
\endsubhead
Suppose that $a=\Cal Ff$ and $b=\Cal Fg$ for 
$f,g\in \Cal O_C(\Bbb R^2)$. Then we have in the weak sense (as 
distributions)
$$\align
&\Cal F\i((\Cal Ff)*_\hbar(\Cal Fg))(x) =\\
&= \tfrac1{(2\pi)^2}\int_{\Bbb R^4}(\Cal Ff)(y-z)(\Cal Fg)(z)
     e^{i(\langle x,y\rangle +(z_1y_2-z_1z_2)\hbar)}dy\,dz \\
&= \tfrac1{(2\pi)^2}\int_{\Bbb R^4}(\Cal Ff)(y)(\Cal Fg)(z)
     e^{i(\langle x,y+z\rangle +z_1y_2\hbar)}dy\,dz \\
&= \tfrac1{(2\pi)^2}\int_{\Bbb R^4}f(u)g(v)
     \left(\int_{\Bbb R^4} 
     e^{i(\langle y,x-u\rangle+\langle z,u-v\rangle 
     +(z_1y_2-z_1z_2)\hbar)}dy\,dz \right)du\,dv\\
\endalign$$
Let us now use $\Cal F=\Cal F_1\o\Cal F_2$, the composition of the 
two one dimensional Fourier transforms in both variables separately, 
and recall that the integals above are weak, are in 
$\Cal O_M(\Bbb R^2)\subset \Cal S'(\Bbb R^2)$, so they make sense 
only when applied to test functions in $\Cal S$. Then the last but 
one expression becomes
$$\align
&= \tfrac1{(2\pi)^2}\int_{\Bbb R^4}(\Cal F_1\Cal F_2f)(y_1,y_2)
     (\Cal F_1\Cal F_2g)(z_1,z_2)\\
&\qquad\qquad\qquad e^{i(x_1y_1+x_2y_2+x_1z_1+x_2z_2+z_1y_2\hbar)}
     dy_1\,dy_2\,dz_1\,dz_2 \\
&= \int_{\Bbb R^2}
     \tfrac1{2\pi}\int_\Bbb R(\Cal F_1\Cal F_2f)(y_1,y_2)e^{ix_1y_1}dy_1\\
&\qquad\qquad\qquad
     \tfrac1{2\pi}\int_\Bbb R(\Cal F_2\Cal F_1g)(z_1,z_2)e^{ix_2z_2}dz_2
     e^{i(x_2y_2+x_1z_1+z_1y_2\hbar)}
     dy_2\,dz_1 \\
&= \int_{\Bbb R^2}(\Cal F_2f)(x_1,y_2)(\Cal F_1g)(z_1,x_2)
     e^{i(x_2y_2+x_1z_1+z_1y_2\hbar)}
     dy_2\,dz_1 \\
&= \sum_{k=0}^\infty\frac{(i\hbar)^k}{k!}
     \int_{\Bbb R}y_2^k(\Cal F_2f)(x_1,y_2)e^{ix_2y_2}dy_2
     \int_\Bbb R z_1^k(\Cal F_1g)(z_1,x_2)e^{ix_1z_2}dz_2\\
&= 4\pi^2\sum_{k=0}^\infty\frac{(-i\hbar)^k}{k!}
     \partial_2^kf(x_1,x_2)
     \partial_1^kg(x_1,x_2),\\
\endalign$$
where we used $i\partial_xf(x)=\Cal F_y\i(y(\Cal Ff)(y))(x)$. 
The last expression is half of the Moyal star product, represented by a 
convergent integral. Obviously the series can only be interpreted as 
a formal power series in $\hbar$. But note that the divergence 
appears only after the interchange of the sum with the integral; 
before the expressions are bounded bilinear in $f$ and $g$, and even 
smooth in $\hbar$.

Also one should compare this result with 
the treatment of the Weyl calculus in \cit!{21},~III, 18.5.

\subhead\nmb.{3.6}. The Fourier transform of the other twisted 
convolution \endsubhead
Let us apply the other twisted convolution to $a=\Cal Ff, 
b=\Cal Fg\in \Cal O_M'(\Bbb R^2)$ for 
$f,g\in \Cal O_C(\Bbb R^2)$: 
$$\align
&\Cal F\i((\Cal Ff)\hat*_\hbar(\Cal Fg))(x) =\\
&= \tfrac1{(2\pi)^2}\int_{\Bbb R^2}((\Cal Ff)\hat*_\hbar(\Cal Fg))(y)
     e^{i\langle x,y\rangle}dy \\
&= \tfrac1{(2\pi)^2}\int_{\Bbb R^4}(\Cal Ff)(y-z)(\Cal Fg)(z)
     e^{i(\langle x,y\rangle -\frac\hbar 2\om(y,z))}dy\,dz \\
&= \tfrac1{(2\pi)^2}\int_{\Bbb R^4}(\Cal Ff)(y)(\Cal Fg)(z)
     e^{i(\langle x,y+z\rangle -\frac\hbar 2\om(y,z))}dy\,dz \\
&= \tfrac1{(2\pi)^2}\int_{\Bbb R^4}(\Cal Ff)(y)e^{i\langle x,y\rangle}
     (\Cal Fg)(z)e^{i\langle x,z\rangle}
     \left(\sum_{k=0}^\infty \frac{(-i\hbar)^k}{2^k\,k!} 
     (y_1x_2-y_2z_1)^k \right)dy\,dz\\ 
\endalign$$
Let us now use 
$(i\partial_1)^m(i\partial_2)^nf(x)
     =\Cal F_y\i(y_1^my_2^n(\Cal Ff)(y))(x)$, 
which also holds in the weak sense for tempered distributions.
Then we may continue to compute in the weak sense of distributions:
$$\align
&= \tfrac1{(2\pi)^2}\sum_{k=0}^\infty \frac{(-i\hbar)^k}{2^k\,k!}
     \int_{\Bbb R^4}(\Cal Ff)(y)e^{i\langle x,y\rangle}
     (\Cal Fg)(z)e^{i\langle x,z\rangle}
     (y_1x_2-y_2z_1)^k dy_1\,dy_2\,dz_1\,dz_2\\ 
&= (2\pi)^2\left.\sum_{k=0}^\infty \frac{(-i\hbar)^k}{2^k\,k!}
     (\partial_{y_2}\,\partial_{z_1}
     -\partial_{y_1}\,\partial_{z_2})^k(f(y_1,y_2)g(z_2,z_2))\right|
     \Sb y_1=z_1=x_1 \\ y_2=z_2=x_2 \endSb.
\endalign$$
This is now really the 
Moyal star product, expressed as a sum of bidifferential operators.

\subhead\nmb.{3.7}. Convolution algebras on the Heisenberg group 
\endsubhead 
Let us consider the Heisenberg group in the following form: 
$\operatorname{He}_2^\hbar = \Bbb R^2 \x S^1$ with multiplication
$$(x,\al).(y,\be)=(x+y,\al\be e^{\frac{i\hbar}2\om(x,y)})
     =(x+y,\al\be e^{\frac{i\hbar}2(x_1y_2-x_2y_1)}).$$
Let us consider the bounded linear mapping between the spaces of 
speedily decreasing distributions
$$\gather
\tilde{\quad}: 
\Cal O_M'(\Bbb R^2)\to \Cal O_M'(\operatorname{He}_2^\hbar),\quad
\tilde a(x_1,x_2,\al)=a(x_1,x_2)\al
\endgather$$
Since the Haar measure on $\operatorname{He}_2^\hbar$ is just the 
usual measure $dx_1\wedge dx_2\wedge d\al$, where we choose
$\int_{S_1}d\al = 1$, we can then compute the convolution as a 
weak integral (in the sense of tempered distributions):
$$\align
(\tilde a * \tilde b)(x,\al) &= \int_{\operatorname{He}_2^\hbar} 
\tilde a(u,\be)\tilde b((x,\al)(u,\be)\i)du d\be \\
&= \int_{\operatorname{He}_2^\hbar} 
     \tilde a(u,\be)\tilde b(x-u,\al\be\i 
     e^{\frac{i\hbar}2(u_1x_2-u_2x_1)})du_1\,du_2\,d\be\\
&= \int_{\operatorname{He}_2^\hbar} 
     a(u_1,u_2)\be b(x_1-u_1,x_2-u_2)\al\be\i 
     e^{\frac{i\hbar}2(u_1x_2-u_2x_1)}du_1\,du_2\,d\be\\
&= \int_{\Bbb R^2} 
     a(u_1,u_2)b(x_1-u_1,x_2-u_2) 
     e^{\frac{i\hbar}2(u_1x_2-u_2x_1)}du_1\,du_2\;\al\\
&= \widetilde{(a \hat *_\hbar b)} (x,\al)
\endalign$$
The groups $\operatorname{He}_2^\hbar$ are all isomorphic for 
$\hbar\ne0$, an isomorphism 
$\operatorname{He}_2^\hbar\to \operatorname{He}_2^1$ is given by 
$(x_1,x_2,\al)\mapsto (\hbar x_1,\hbar x_2,\al)$. 
Thus all the algebras $(\Cal O_M',\hat *_\hbar)$ are isomorphic for 
$\hbar\ne 0$, in strong contrast to the behaviour of the 
subalgebras $T^2_{e^{i\hbar}}$, the noncommutative tori.

\subhead\nmb.{3.8}. Derivations \endsubhead
Let us determine all derivations of the smooth Heisenberg plane. We 
use the form $(\Cal O_M',\hat *_\hbar)$ from \nmb!{3.3.5}, and 
we start with the inner derivations.
We have for $a,b\in\Cal O_M'(\Bbb R^2)$ in the weak sense
$$\align
\ad(a)b 
&=\int_{\Bbb R^2} a(x-y)b(y) e^{-\frac{i\hbar}2\om(x,y)}\,dy
-\int_{\Bbb R^2} b(x-y)a(y) e^{-\frac{i\hbar}2\om(x,y)}\,dy\tag{\nmb:{1}}\\
&=\int_{\Bbb R^2} a(x-y)b(y) \left(e^{-\frac{i\hbar}2\om(x,y)}
     +e^{\frac{i\hbar}2\om(x,y)}\right)\,dy\\
&=\int_{\Bbb R^2} a(x-y)b(y) 
     2\cos\left(\tfrac{\hbar}2\om(x,y)\right)\,dy
\endalign$$

\proclaim{Proposition} Every bounded derivation of 
$(\Cal O_M'(\Bbb R^2),*_\hbar)$ is inner, if $\hbar\ne0$.
\endproclaim

\demo{Proof}
Let us note first that $Q$ and $P$ are elements of $\Cal O_M'$, 
namely we have  
$$
Q=\int\de'_t(0)\de_s(0)e^{itQ}e^{isP}\,dt\,ds,\quad\text{etc.}
$$
Let 
$D:(\Cal O_M'(\Bbb R^2),*_\hbar)\to (\Cal O_M'(\Bbb R^2),*_\hbar)$ be 
a bounded derivation. 
Then we let 
$$\gather
D(Q) = \int a_Q(t,s)e^{itQ}e^{isP}dt\,ds,\quad
D(P) = \int a_P(t,s)e^{itQ}e^{isP}dt\,ds
\endgather$$
We want to find a distribution $b\in\Cal O_M'$ such that 
$B=\int b(t,s)e^{itQ}e^{isP}dt\,ds$ satisfies $D(Q)=[B,Q]$ and 
$D(P)=[B,P]$.
We have (using formulas from the analytic proof of lemma \nmb!{3.2}) 
for $f\in \Cal S(\Bbb R)$
$$\align
([e^{isP},Q]f)(u) &= (zf(z))|_{z=u+s\hbar} - uf(u+s\hbar)
= (s\hbar e^{isP}f)(u)\\
[B,Q] &= \int b(t,s)e^{itQ}[e^{isP},Q]\,dt\,ds
     = \int b(t,s)s\hbar e^{itQ}e^{isP}\,dt\,ds,
\endalign$$
and similarly
$$\align
([e^{itQ},P]f)(u) &= e^{itu}\tfrac{\hbar}i\partial_uf(u) - 
\tfrac{\hbar}i\partial_u(e^{itu}f(u))
= -(t\hbar e^{itQ}f)(u)\\
[B,P] &= \int b(t,s)[e^{itQ},P]e^{isP}\,dt\,ds
     = -\int b(t,s)t\hbar e^{itQ}e^{isP}\,dt\,ds,
\endalign$$
so that we have to solve 
$$
b(t,s)s\hbar = a_Q(t,s), \quad -b(t,s)t\hbar= a_P(t,s).
$$
Applying the Fourier transform we have to find $\hat b\in \Cal O_C$ 
which satisfies
$$\gather
i\hbar \partial_s \hat b(t,s) = \hat a_Q(s,t), \quad
i\hbar \partial_t \hat b(t,s) = -\hat a_P(s,t),\\
d \hat b = \tfrac1{i\hbar}(\hat a_Q ds - \hat a_P dt).
\endgather$$
This can be solved by the Lemma of Poincar\'e in $\Cal O_C(\Bbb R^2)$ 
if and only if 
$$d(\hat a_Q ds - \hat a_P dt)
= (\partial_t\hat a_Q + \partial_s\hat a_P) dt\wedge ds = 0.$$
But this is the case since we have in turn, using the results from 
above,
$$\gather
D(Q)P+QD(P)-D(P)Q-PD(Q) = D([Q,P]) = D(i\hbar.1) = 0\\
[D(Q),P] = [D(P),Q]\\
-a_Q(t,s)t\hbar = a_P(t,s)s\hbar\\
i\hbar\partial_t \hat a_Q + i\hbar\partial_s \hat a_P =0,
\endgather$$
as required. The lemma of Poincar\'e has the form: $d\ph=0$ implies 
$\ph=d\ps$ where $\ps(x)=\int_0^1\sum_i\ph_i(tx)x_idt$. Thus 
$\ph_i\in \Cal O_C(\Bbb R^2)$ implies $\ps\in \Cal O_C(\Bbb R^2)$ by a 
simple estimation. 
Thus on $Q$ and $P$ the bounded derivation $D$ agrees 
with an inner derivation.

It remains to show that a bounded derivation $D$ which vanishes on 
$Q$ and on $P$ must vanish on $\Cal O_M'$.
For that we note the following facts:

The curve $t\mapsto e^{itQ}$ is a smooth 1-parameter group of isomorphisms of 
$\Cal S(\Bbb R)$ with infinitesimal generator $iQ$, and it is the 
unique 1-parameter group with this generator, since for any other 
$C(t)$ we have 
$\partial_t(e^{itQ})C(-t)=e^{itQ}iQC(t)-e^{itQ}iQC(t)=0$, so that 
$e^{itQ}C(-t)$ is the constant $Id$.

Consider the semidirect product 
$(\Cal O_M'(\Bbb R^2),*_\hbar)\ltimes \Cal O_M'(\Bbb R^2)\ep$
where $\ep$ is in the center and $\ep^2=0$, with the multiplication 
$(a+b\ep).(a'+b'\ep)=aa'+(ab'+ba')\ep$.
Obviously $D$ is a derivation if and only if $a\mapsto a+D(a)\ep$ is 
a homomorphism of algebras. 

Thus $t\mapsto e^{itQ}+D(e^{itQ})\ep$ is a smooth 1-parameter group in 
the semidirect product with infinitesimal generator 
$iQ+D(iQ)\ep=iQ+0$ and with second 1-parameter group $e^{itQ}+0$, 
thus $D(e^{itQ})=0$ for all $t$.

Similarly $D(e^{isP})=0$ for all $s$.
Thus $D$ vanishes on $e^{itQ}e^{isP}$ 
for each $t$ and $s$. And if 
$a\in C^\infty_c(\Bbb R^2)\subset \Cal O_M'(\Bbb R^2)$, then 
$$
D\left( \int a(t,s)e^{itQ}e^{isP}dt\,ds\right) = 
\int a(t,s)D\left(e^{itQ}e^{isP}\right)dt\,ds = 0,
$$
since Riemann sums converge Mackey to the integral.
Finally one should note that $C^\infty_c$ is dense in $\Cal O_M'$, so the 
result follows.
\qed\enddemo

\proclaim{\nmb.{3.9}. Conjecture} 
The smooth Heisenberg plane $\Cal O_M'(\Bbb R^{2n})$ is a smooth 
$*$-algebra with derivation space the space of all bounded 
derivations in the given topology, and a suitable state space. \endproclaim

In fact we think that the topology described in \nmb!{1.1.3} is 
the one of $\Cal O_M'(\Bbb R^{2n})\cong \Cal O_C(\Bbb R^{2n})$.
One has to show that each state is a bounded linear functional, and 
that we are able to find enough states and derivations in order to 
generate the topology described in \nmb!{3.1}.

\heading\totoc\nmb0{4}. APPENDIX: Calculus in infinite dimensions and 
convenient vector spaces \endheading

\subheading{\nmb.{4.1}} 
The notion of convenient vector spaces arose in the quest for the 
right setting for differential calculus in infinite dimensions:
The traditional approach to differential calculus works 
well for Banach spaces, but for more general locally convex spaces 
there are difficulties. The main one is that 
the composition of linear mappings stops to be jointly continuous at 
the level of Banach spaces, for any compatible topology, so that 
even the chain rule is not valid without further assumptions.
In addition to their importance for differential calculus convenient 
vector spaces together with bounded linear mappings and the 
appropriate tensor product form a monoidally closed category, the 
only useful one which functional analysis offers beyond Banach spaces.

In this section we sketch the basic definitions and the most 
important results concerning calculus for convenient vector spaces.
All locally convex spaces will be 
assumed to be Hausdorff.
Proofs for the results sketched here can be found in 
\cit!{12} (sauf for \nmb!{4.8} which was proved in \cit!{5}). 
A complete coverage is in the book 
\cit!{20}; \cit!{5} contains an overview and a presentation of 
non-commutative geometry based on convenient vector spaces. 

\subheading{\nmb.{4.2}. Smooth curves} Let $E$ be a locally 
convex vector space. A curve $c:\Bbb R\to E$ is called {\it smooth} 
or $C^\infty$ if all derivatives exist (and are continuous) - this is 
a concept without problems. Let $C^\infty(\Bbb R,E)$ be the space of 
smooth curves. It can be shown that the set $C^\infty(\Bbb R,E)$ does depend 
on the locally convex topology of $E$ only through its
underlying bornology (system of bounded sets). 
\comment
The final topologies 
with respect to the following sets of mappings into $E$ coincide:
\roster
\item"(\nmb:{1})" $C^\infty(\Bbb R,E)$.
\item"(\nmb:{2})" Lipschitz curves (so that $\{\frac{c(t)-c(s)}{t-s}:t\neq s\}$ 
     is bounded in $E$). 
\item"(\nmb:{3})" $\{E_B\to E: B\text{ bounded absolutely convex in }E\}$, where 
     $E_B$ is the linear span of $B$ equipped with the Minkowski 
     functional $p_B(x):= \inf\{\la>0:x\in\la B\}$.
\item"(\nmb:{4})" Mackey-convergent sequences $x_n\to x$ (there exists a sequence 
     $0<\la_n\nearrow\infty$ with $\la_n(x_n-x)$ bounded).
\endroster
This topology is called the $c^\infty$-topology on $E$ and we write 
$c^\infty E$ for the resulting topological space. In general (on the 
space $\Cal D$ of test functions for example) it is finer than the 
given locally convex topology; it is not a vector space topology, 
since addition is no longer jointly continuous. The finest among all 
locally convex topologies on $E$ which are coarser than the 
$c^\infty$-topology is the bornologification of the given locally 
convex topology. If $E$ is a Fr\'echet space, then $c^\infty E = E$.
\endcomment

\subheading{\nmb.{4.3}. Convenient vector spaces} Let $E$ be a 
locally convex vector space. $E$ is said to be a {\it convenient 
vector space} if one of the following equivalent
conditions is satisfied (called $c^\infty$-completeness):
\roster
\item"(\nmb:{1})" Any Mackey-Cauchy-sequence (so that there are
     scalars $\la_{n,m}\to\infty$ such that 
		 $\{\la_{n,m}(x_n-x_m):n,m\in\Bbb N\}$ is bounded) 
     converges. 
\item"(\nmb:{2})" If $B$ is bounded closed and absolutely convex, 
     then the linear span $E_B$ of $B$ is a 
     Banach space with respect to the Minkowski functional
		 $p_B(x):= \inf\{\la>0:x\in\la B\}$.
\item"(\nmb:{3})" Any Lipschitz curve (so that $\{\frac{c(t)-c(s)}{t-s}:t\neq s\}$ 
     is bounded) in $E$ is locally Riemann integrable.
\item"(\nmb:{4})" For any $c_1\in C^\infty(\Bbb R,E)$ there is 
     $c_2\in C^\infty(\Bbb R,E)$ with $c_1=c_2'$ (existence of 
     antiderivative).
\item"(\nmb:{5})" If $f:\Bbb R\to E$ is scalarwise $\Lip^k$, then $f$ is 
     $\Lip^k$, for $k>1$.
\item"(\nmb:{6})" If $f:\Bbb R\to E$ is scalarwise $C^\infty$ then $f$ is 
     differentiable at 0.
\item"(\nmb:{7})" If $f:\Bbb R\to E$ is scalarwise $C^\infty$ then $f$ is 
     $C^\infty$.
\endroster
Here a mapping $f:\Bbb R\to E$ is called $\Lip^k$ if all partial 
derivatives up to order $k$ exist and are Lipschitz, locally on 
$\Bbb R$. To be scalarwise $C^\infty$ means for a curve $f$ 
that $\la\o f$ is $C^\infty$  
for all continuous (equivalently: all bounded) linear functionals $\la$ on $E$.
Obviously $c^\infty$-completeness is weaker than 
sequential completeness, so any sequentially complete locally convex
vector space is convenient.
 From \nmb!{4.2.4} one easily sees that (sequentially) closed linear
subspaces of convenient vector spaces are again convenient. We 
always assume that a convenient vector space is equipped with its 
bornological topology. All spaces which a working mathematician needs 
in daily life are convenient.
For any locally convex space $E$ there is a convenient vector space
$\tilde E$ called the completion of $E$, and a bornological embedding 
$i:E\to \tilde E$, which is characterized by the
property that any bounded linear map from $E$ into an arbitrary 
convenient vector space extends to $\tilde E$.

\subheading{\nmb.{4.4}. Smooth mappings} Let $E$ and $F$ be locally 
convex vector spaces.
A mapping $f:E\to F$ is called {\it smooth} or 
$C^\infty$, if $f\o c\in C^\infty(\Bbb R,F)$ for all 
$c\in C^\infty(\Bbb R,E)$; so 
$f_*: C^\infty(\Bbb R,E)\to C^\infty(\Bbb R,F)$ makes sense.
Let $C^\infty(E,F)$ denote the space of all smooth mappings from $E$ 
to $F$. For $E$ and $F$ finite dimensional (or even Fr\'echet spaces) 
this gives the usual notion of smooth mappings  (Already for $E=\Bbb R^2$ 
this is a non-trivial statement). 
Multilinear mappings are smooth if and 
only if they are bounded. We denote by $L(E,F)$ the space 
of all bounded linear mappings from $E$ to $F$. 

\subheading{\nmb.{4.5}. Differential calculus} We equip the 
space $C^\infty(\Bbb R,E)$ with the bornologification of the topology 
of uniform convergence on compact sets, in all derivatives 
separately. Then we equip the space $C^\infty(E,F)$ with the 
bornologification of the initial topology with respect to all 
mappings $c^*:C^\infty(E,F)\to C^\infty(\Bbb R,F)$, $c^*(f):=f\o c$, 
for all $c\in C^\infty(\Bbb R,E)$.
We have the following results:
{\sl
\roster
\item"(\nmb:{1})" If $F$ is convenient, then also $C^\infty(E,F)$ is convenient, 
     for any $E$. The space $L(E,F)$ is a closed linear subspace of 
     $C^\infty(E,F)$, so it is convenient also.
\item"(\nmb:{2})" The smooth uniform boundedness principle: 
     If $E$ is convenient, then a curve $c:\Bbb R\to L(E,F)$ is 
     smooth if and only if $t\mapsto c(t)(x)$ is a smooth curve in $F$ 
     for all $x\in E$.
\item"(\nmb:{3})" The category of convenient vector spaces and smooth mappings is 
     cartesian closed. So we have a natural bijection 
$$C^\infty(E\x F,G)\cong C^\infty(E,C^\infty(F,G)),$$
     which is even a homeomorphism. 
     Note that this result, for $E=\Bbb R$, is the prime assumption of 
     variational calculus.
     As a consequence evaluation mappings, insertion mappings, and 
     composition are smooth.
\item"(\nmb:{4})" The differential $d: C^\infty(E,F)\to C^\infty(E,L(E,F))$,
     given by $df(x)v:=\lim_{t\to0}\frac1t(f(x+tv)-f(x))$,
     exists and is linear and bounded (smooth). Also the chain rule holds: 
     $d(f\o g)(x)v = df(g(x))dg(x)v$.
\endroster
}

\subheading{\nmb.{4.6} }
The category of convenient vector spaces and bounded linear maps
is complete and cocomplete, so all categorical limits and colimits
can be formed. In particular we can form products and direct sums of
convenient vector spaces.

For convenient vector spaces $E_1$,\dots ,$E_n$, and $F$ we can now
consider the space of all bounded $n$-linear maps, $L(E_1,\dots ,E_n;F)$,
which is a closed linear subspace of $C^\infty(\prod _{i=1}^nE_i,F)$
and thus again convenient. It can be shown that multilinear maps are
bounded if and only if they are partially bounded, i\.e\. bounded in 
each coordinate and that there is a natural isomorphism (of 
convenient vector spaces) $L(E_1,\dots ,E_n;F)\cong
L(E_1,\dots ,E_k;L(E_{k+1},\dots ,E_n;F))$

\proclaim{\nmb.{4.7}. Result}
On the category of convenient vector spaces there is a unique tensor product
$\tilde \otimes$ which makes the category symmetric monoidally 
closed, i\.e\. there are natural isomorphisms of convenient vector 
spaces $L(E_1;L(E_2,E_3))\cong
L(E_1\tilde \otimes E_2,E_3)$, 
$E_1\tilde \otimes E_2\cong E_2\tilde \otimes E_1$,
$E_1\tilde \otimes (E_2\tilde \otimes E_3)\cong (E_1\tilde \otimes
E_2)\tilde \otimes E_3$ and $E\tilde \otimes \Bbb R\cong E$.
\endproclaim

\proclaim{\nmb.{4.8}. Result} {\rm \cit!{5}, 2.7}.
Let $A$ be a convenient algebra, $M$ a convenient right $A$-module
and $N$ a convenient left $A$-module. This means that all structure 
mappings are bounded bilinear.
\roster
\item"(\nmb:{1})" There is a convenient vector space 
     $M\tilde \otimes _AN$ and a bounded bilinear map 
     $b:M\times N\to M\tilde \otimes _AN$, 
     $(m,n)\mapsto m\otimes _An$ such that 
     $b(ma,n)=b(m,an)$ for all $a\in A$, $m\in M$ and $n\in N$ which has
     the following universal property: If $E$ is a convenient vector space
     and $f:M\times N\to E$ is a bounded bilinear map such that 
     $f(ma,n)=f(m,an)$ then there is a unique bounded linear map
     $\tilde f:M\tilde \otimes _AN\to E$ with $\tilde f\o b=f$.
\item"(\nmb:{2})" Let $L^A(M,N;E)$ denote the space of all bilinear 
     bounded maps $f:M\times N\to E$ having the above property, which 
     is a closed linear subspace of $L(M,N;E)$. Then we have an 
     isomorphism of convenient vector spaces 
     $L^A(M,N;E)\cong L(M\tilde \otimes _AN,E)$. 
\item"(\nmb:{3})" If $B$ is another convenient algebra such that $N$ is a 
     convenient right $B$-module and such that the actions of $A$
     and $B$ on $N$ commute, then $M\tilde \otimes _AN$ is in a canonical
     way a convenient right $B$-module.
\item"(\nmb:{4})" If in addition $P$ is a convenient left $B$-module 
     then there is a natural isomorphism of convenient vector spaces
$$M\tilde \otimes _A(N\tilde \otimes _BP)\cong (M\tilde \otimes 
     _AN)\tilde \otimes _BP$$
\endroster
\endproclaim

\enddocument